\documentclass[12pt]{amsart}
\usepackage{amssymb}
\author{Vitaly Bergelson and Alexander Gorodnik}
\newtheorem{thm}{Theorem}
\newtheorem{lem}[thm]{Lemma}
\newtheorem{pro}[thm]{Proposition}
\newtheorem{cor}[thm]{Corollary}

\newtheorem{rmk}[thm]{Remark}

\newtheorem{que}[thm]{Question}
\newtheorem{exa}[thm]{Example}

\newcommand\Aut{\operatorname{Aut}}

\newcommand{\ignore} [1] {}
\numberwithin{equation}{section}
\numberwithin{thm}{section}

\title{Ergodicity and mixing of noncommuting epimorphisms}
\thanks{The first author was partially supported by
  NSF grant DMS-0345350.}
\thanks{ The second author was partially supported by NSF grant DMS-0400631.}
\address{Department of Mathematics, The Ohio State University,
  Columbus, OH 43210, USA}
\email{vitaly@math.ohio-state.edu}   
\address{Mathematics 253-37, CalTech, Pasadena, CA 91125, USA}
\email{gorodnik@caltech.edu}

\begin{document}
\begin{abstract}
We study mixing properties of epimorphisms of a compact connected
finite-dimensional abelian group $X$. In particular, we show that a
set $F$, $|F|>\dim X$, of epimorphisms of $X$ is mixing iff every
subset of $F$ of
cardinality $(\dim X)+1$ is mixing. We also construct examples of
free nonabelian groups of automorphisms of tori which are mixing,
but not mixing of order 3, and show that, under some irreducibility
assumptions, ergodic groups of automorphisms contain mixing subgroups
and free nonabelian mixing subsemigroups.
\end{abstract}

\maketitle

\section{Introduction}

\subsection{Mixing sets}
Let $X$ be a compact abelian group, $\mathcal{B}$ the completion of
the Borel $\sigma$-algebra of $X$, and $m$
the normalized Haar measure on $X$. A finite set $F$, $|F|>1$, of
epimorphisms (i.e., continuous surjective self-homomorphisms) 
of $X$ is
called {\it mixing} if for any collection of measurable sets $B_\gamma\in\mathcal{B}$, $\gamma\in F$,
$$
m\left(\bigcap_{\gamma\in F} \gamma^{-n}(B_\gamma)\right)\to \prod_{\gamma\in F} m(B_\gamma)\quad\hbox{as $n\to\infty$}.
$$
Such set is sometimes also called mixing shape.
It is clear that if $F$ is mixing, then every subset of $F$ is mixing as well.
However, in  general, the assumption that all proper subsets of $F$
are mixing does not imply that $F$ is mixing.
For example, it was shown by F.~Ledrappier that there exist 
commuting automorphisms $\gamma_1$ and $\gamma_2$ of a compact totally
disconnected abelian group such that the sets $\{id,\gamma_1\}$
$\{id,\gamma_2\}$,  $\{\gamma_1,\gamma_2\}$ are mixing, but the set $\{id,\gamma_1,\gamma_2\}$
is not mixing (see \cite{led} and \cite[Chapter~VIII]{sch}).
Also, if one does not assume commutativity, similar examples exist
  for connected groups as well (see Corollary \ref{th:ex} below).

K.~Schimdt has shown that when the group $X$ is
connected and the epimorphisms which form the set $F$ commute,
the situation is quite different  (see \cite{sch0}):

\begin{thm}[Schmidt]\label{th:comm}
Let $X$ be a compact connected abelian group and $F$ a finite set of commuting epimorphisms of $X$. 
Then the set $F$ is mixing iff every subset of $F$ of cardinality $2$ is mixing.
\end{thm}

In this paper, we prove a noncommutative analog of Theorem \ref{th:comm}:

\begin{thm}\label{th_d+1}
Let $X$ be a compact connected abelian group such that $\dim X=d<\infty$
and $F$ a finite set of epimorphisms of $X$ with $|F|>\dim X$.
Then the set $F$ is mixing iff every subset $E$ of $F$ with $|E|=d+1$ is mixing.
\end{thm}

Theorem \ref{th_d+1} and Theorem \ref{th:comm}
(in the finite-dimensional case) follow from Theorem \ref{th:general} below.
We also show that the bound $d+1$ in Theorem \ref{th_d+1} is sharp (see
Corollary \ref{th:ex} below).

\subsection{Mixing sets and spectrum}\label{sec:spectrum}

Let $X$ be a compact connected abelian group with $\dim X=d<\infty$.
We denote by $\hat X$ the character group of $X$. 
Under the above assumptions, $\hat X$ is a discrete abelian torsion
free group of rank $d$. Hence, we may assume that
$$
\mathbb{Z}^d\subset \hat X\subset \mathbb{Q}^d.
$$
(Conversely, any abelian group $A$ such that $\mathbb{Z}^d\subset A\subset \mathbb{Q}^d$
corresponds to a compact connected abelian group of dimension $d$.)

Any continuous endomorphism $T$ of $X$ defines an endomorphism $\hat
T$ of $\hat X$ that extends to a linear map of $\mathbb{Q}^d$.
Note that $T$ is surjective iff $\hat T$ is nondegenerate
(i.e., $\det \hat T\ne 0$).

We establish the following criterion  for mixing in terms of eigenvectors
of the corresponding linear maps of $\mathbb{Q}^d$.

\begin{thm}\label{th:general}
A set $\{T_1,\ldots,T_s\}$ of epimorphisms of $X$ is mixing iff
for every $l\ge 1$, every subset $\{k_1,\ldots,k_r\}\subset \{1,\ldots,s\}$,
and every $\lambda\in\mathbb{C}$,
there are no $\lambda$-eigenvectors of $\hat T_{k_1}^l,\ldots, \hat T_{k_r}^l$
that are linearly dependent over $\mathbb{Q}$.
\end{thm}

\begin{rmk}{\rm
It follows from the proof that in Theorem \ref{th:general} one can
replace ``for every $l\ge 1$'' by ``for every $l\ge 1$ such that
$\phi(l)\le (\dim X)^2$'', where $\phi$ denotes the Euler's totient
function. Moreover, this estimate is sharp (see Example \ref{ex:last2} below).
}
\end{rmk}

We state some corollaries of Theorem \ref{th:general}.
Note that Corollary \ref{cor:com} is just another formulation of Theorem \ref{th:comm}
in the finite-dimensional case, and Corollary \ref{th:cond_gen} implies Theorem \ref{th_d+1}.

\begin{cor}\label{cor:com}
For commuting epimorphisms $T_1,\ldots, T_s$ of $X$, the following statements are
equivalent:
\begin{enumerate}
\item[(a)] The set $\{T_1,\ldots, T_s\}$ is mixing.
\item[(b)] For every $i\ne j$, the set $\{T_i,T_j\}$ is mixing.
\item[(c)] For every $i\ne j$, the linear map $\hat T_i^{-1}\hat T_j$
does not have roots of unity as eigenvalues.
\end{enumerate}
\end{cor}

For two (not necessarily commuting) epimorphisms, we have the
following criterion for mixing:

\begin{cor}\label{th:two}
The set of epimorphisms $\{T_1,T_2\}$ of $X$ is mixing iff
there is no closed subgroup $Y\ne X$ such that 
for some $l\ge 1$,  $Y$ is $\{T_1^l,T_2^l\}$-invariant and $T_1^l=T_2^l$ on $X/Y$.
\end{cor}

Corollary \ref{th:two} may fail if the group $X$ is
disconnected or infinite-dimensional (see Example \ref{ex:two} below).


Denote by $\hbox{Spec}(T)$ the set of eigenvalues of $\hat T$. 
The following corollary of Theorem \ref{th:general} characterizes
mixing in terms of spectrum:

\begin{cor}\label{th:cond_gen}
Let $T_1,\ldots, T_s$ be epimorphisms of $X$.
\begin{enumerate}
\item[(a)] If for every $l\ge 1$ and $i,j=1,\ldots,s$, $i\ne j$,
\begin{align*}
\hbox{\rm Spec}(T_i^l)\cap\hbox{\rm Spec}(T_j^l)=\emptyset,
\end{align*}
then $\{T_1,\ldots, T_s\}$ is mixing.
\item[(b)] If for some $l\ge 1$ and $S\subset \{1,\ldots,s\}$ such that  $|S|>d$,
$$
\bigcap_{i\in S} \hbox{\rm Spec}(T_i^l) \ne \emptyset,
$$
then $\{T_1,\ldots, T_s\}$ is not mixing.
\item[(c)]  If for every $l\ge 1$ and $S\subset \{1,\ldots,s\}$ such that  $|S|>d$,
$$
\bigcap_{i\in S} \hbox{\rm Spec}(T_i^l) = \emptyset,
$$
then 
$\{T_1,\ldots, T_s\}$ is mixing iff every subset of cardinality $d$ is mixing.
\end{enumerate}
\end{cor}

\begin{rmk}{\rm
In Corollary \ref{th:cond_gen}, one can
replace ``for every $l\ge 1$'' by ``for every $l\ge 1$ such that
$\phi(l)\le (\dim X)^2$''. Moreover, this estimate is sharp (see Example \ref{ex:last2} below).
}
\end{rmk}

Corollary \ref{th:cond_gen}(a) shows that, if
epimorphisms $T_1,\ldots,T_s$ are ``spectrally independent'',
then for every 
$B_1,\ldots, B_s\in\mathcal{B}$,
$$
\lim_{n\to\infty} m(T_1^{-n} B_1\cap\cdots \cap T_s^{-n} B_s)=m(B_1)\cdots m(B_s).
$$
Although this limit does not exist in general (consider, for example, $T_1=id$ and $T_2=-id$),
the proof of Theorem \ref{th:general} implies the following corollary.

\begin{cor}\label{thm:limit}
For any finite set $\{T_1,\ldots, T_s\}$ of epimorphisms of $X$, there exists $l\ge 1$
such that for every $k\in\mathbb{Z}/l\mathbb{Z}$ and $f_1,\ldots, f_s\in L^\infty(X)$, the limit
\begin{equation}\label{eq:limit}
\lim_{{\tiny\begin{matrix} n\to\infty \\ n=k(\hbox{\rm mod}\, l)\end{matrix}} } \int_X f_1(T_1^{n}x)\cdots f_s(T_s^{n}x)\,dm(x)
\end{equation}
exists.
\end{cor}

\begin{rmk}{\rm
\begin{enumerate}
\item[(i)] It follows from the proof that the integer $l$ appearing in
  Corollary \ref{thm:limit} can be chosen so that $\phi(l)\le (\dim
  X)^2$.  Moreover, this estimate is sharp (see Example \ref{ex:last2} below).

\item[(ii)]
Corollary \ref{thm:limit} is, in general, false  if the group
$X$ is either infinite-dimensional or disconnected (see Example
\ref{ex:limit} below).

\item[(iii)]
Existence of the Ces\`aro limit
$$
\lim_{N-M\to\infty} \frac{1}{N-M}\sum_{n=M+1}^N f_1(T_1^{n}x)\cdots f_s(T_s^{n}x)
$$
in $L^2(X)$ for a certain class
of epimorphisms of a compact abelian group $X$
was proved by D.~Berend in \cite{ber2}. Corollary \ref{thm:limit}
strengthens Berend's result in the case when the  group $X$ is
connected and finite-dimensional.
\end{enumerate}
}
\end{rmk}

We call an automorphism $T$ of $X$ {\it unipotent} if the matrix $\hat T$
is unipotent.

\begin{cor}\label{th:ex}
\begin{enumerate}
\item[(a)]
For every $s=2,\ldots,d+1$ there exists a set $F$ with $|F|=s$
consisting of unipotent
automorphisms of $\mathbb{T}^d$ such that $F$ is not mixing, but every proper subset
of $F$ is mixing.
\item[(b)]
For every $s=2,\ldots, d+1$ there exists a set of mixing
epimorphisms $F$ of $\mathbb{T}^d$ with $|F|=s$
such that $F$ is not mixing, but every proper subset of $F$ is mixing.
\end{enumerate}
\end{cor}


The following corollary relates the notion of ``mixing sets''
(terminology from \cite{sch}) with the notion of ``jointly mixing
automorphisms'' which was introduced in \cite{ber} and used in \cite{bg}.
Epimorphisms $T_1,\ldots, T_{s-1}$ are called {\it jointly
mixing} if the set $\{T_1,\ldots, T_{s-1}, id\}$ is mixing in our terminology.

\begin{cor}\label{p:strong}
The set $\{T_1,\ldots, T_{s-1}, id\}$ of epimorphisms of $X$ is mixing iff
every $T_i$ is mixing and
$\{T_1,\ldots,T_{s-1}\}$ is mixing.
\end{cor}

\subsection{Mixing groups and semigroups}\label{subsec:group}
A semigroup $\Gamma$ of epimorphisms of $X$ is called {\it mixing} if
for every $A,B\in\mathcal{B}$,
$$
m(A\cap \gamma^{-1}B)\to m(A)m(B)
$$
as $\gamma\to\infty$.

A semigroup $\Gamma$ of of epimorphisms of $X$ is {\it mixing of order $s$} if for every $B_1,\ldots, B_s\in\mathcal{B}$,
$$
m(\gamma_1^{-1}B_1\cap(\gamma_2\gamma_1)^{-1}B_2\cap\cdots\cap (\gamma_s\cdots \gamma_1)^{-1}B_s)\to m(B_1)\cdots m(B_s)
$$
as the product $\gamma_j\cdots\gamma_i\to\infty$ for $1<i\le j\le s$.
Note that mixing corresponds to mixing of order 2.

We recall a classical result of Rokhlin (see \cite{roh1}):

\begin{thm}[Rokhlin]\label{th:roh}
If a continuous epimorphism $T$ of a compact abelian group is ergodic,
then it is {\rm mixing of all orders}, that is,
for every $s\ge 1$, $B_1,\ldots, B_s\in\mathcal{B}$, and $n_1,\ldots,n_s\in\mathbb{N}$ such that $|n_i-n_j|\to\infty$
for $i\ne j$,
$$
m(T^{-n_1}B_1\cap\cdots\cap T^{-n_s}B_s)\to m(B_1)\cdots m(B_s).
$$
\end{thm}

This result was extended to finitely generated abelian groups of
automorphisms by K.~Schmidt and T.~Ward in \cite{sw}:

\begin{thm}[Schmidt, Ward]\label{th:scw}
Let $X$ be a compact connected abelian group and $\Gamma\subset\Aut
(X)$, $\Gamma\simeq \mathbb{Z}^n$. Then $\Gamma$ consists of
ergodic automorphisms iff it is mixing of all orders.
\end{thm}

Note that the ergodic properties of the actions in Theorems \ref{th:roh} and
\ref{th:scw} are quite different. The epimorphism $T$ in Theorem \ref{th:roh} has completely
positive entropy (see \cite{roh2}), but the entropy of $\Gamma$-action in
Theorem \ref{th:scw} is
zero if $n>1$ (see \cite[Ch.~V]{sch}). 

While it is true that an arbitrary group $\Gamma$ of automorphisms
is mixing provided that every element of infinite order is
ergodic (see Corollary \ref{c:mix} below), the statement about higher order of mixing fails if $\Gamma$
is not virtually abelian. As an easy corollary of Corollary
\ref{th:cond_gen}(b), we deduce the following result:

\begin{cor}[Bhattacharya]\label{p:bhat}
Let $X$ be a compact connected abelian group with $\dim X=d<\infty$.
Then every subgroup of $\Aut(X)$ which is not virtually abelian is not
mixing of order $d+1$.
\end{cor}

Note that there exist free nonabelian  semigroups of epimorphisms which are
mixing of all orders (see Examples \ref{ex:mix2} and \ref{ex:mix3} below).

Corollary \ref{p:bhat} was first proved by Bhattacharya in \cite{bh}. He also discovered some
interesting rigidity properties of mixing subgroups which are not
virtually abelian. However, it is not obvious whether such subgroups
exist. In this direction, we show:

\begin{thm}
For every $d\ge 2$, $d\ne 3,5,7$, there exists a not virtually abelian
mixing subgroup of $\hbox{\rm Aut}(\mathbb{T}^d)$
which is not mixing of order 3.
\end{thm}

At present, we don't know whether there are such examples for $d=3,5,7$.

Mixing property is much better understood for $\mathbb{Z}^n$-actions
by automorphisms of a compact abelian
group $X$. When $X$ is connected, 2-mixing implies mixing of all orders (see
Theorem \ref{th:scw}). If $X$ is totally disconnected, then 
for every $s\ge 2$, there are examples
that are $s$-mixing but not $(s+1)$-mixing  (see \cite{eiw}).
It is also known that a $\mathbb{Z}^n$-action is $s$-mixing iff every subset of $\mathbb{Z}^n$ of
cardinality $s$ is mixing (see \cite{mas}).

\subsection{Ergodicity and mixing}\label{s:e}

In this subsection we discuss some analogs of Rokhlin's theorem (Theorem \ref{th:roh}) for general groups of
automorphisms. Namely, given a compact abelian group $X$ and a
subgroup $\Gamma$ of $\hbox{Aut}(X)$, we investigate whether ergodicity implies
mixing and mixing of higher orders.
Recall that $\Gamma$ is called {\it
  ergodic} if every measurable $\Gamma$-invariant  subset of $X$ has measure $0$
or $1$. Ergodicity is a weaker notion than mixing. In fact, if $\Gamma$
contains a mixing automorphism, then it is ergodic. 
D.~Berend showed in \cite{ber0} that the converse is also true in the case when $\Gamma$
  is abelian:
\begin{thm}[Berend]\label{th:ber}
Let $X$ be a compact connected finite-dimensional abelian group and
$\Gamma$ an ergodic abelian semigroup of epimorphisms of $X$. Then
$\Gamma$ contains an ergodic   epimorphism.
\end{thm}

Note that by Rokhlin's theorem, an ergodic epimorphism is mixing of all orders.

On the other hand, if $\Gamma$ is not abelian, it may contain no mixing elements (see
  \cite{ber0} or Examples
  \ref{ex:not_irr} and \ref{ex:not_erg} below).
A somewhat stronger version of ergodicity --- ``hereditary
  ergodicity'', which we will presently introduce, is more closely related
  to mixing and will allow us to naturally generalize Berend's theorem.

Let $X$ be a compact abelian group, $Y$ a closed subgroup of $X$, and $\Gamma\subset \hbox{Aut}(X)$.
We define 
$$
\Gamma_Y=\{\gamma\in \Gamma:\, \gamma\cdot Y\subset Y\}.
$$
If $\Gamma_Y$ has finite index in $\Gamma$, we call the
subgroup $Y$ {\it virtually $\Gamma$-invariant}. In the case when $X$
contains no proper closed connected virtually $\Gamma$-invariant subgroups,
we call the group $\Gamma$ {\it strongly irreducible}.
Note that for connected group $X$, strong irreducibility implies ergodicity
(see Proposition \ref{l:irr} below), but the converse is not true (see
Example \ref{ex:not_irr} below).  
We call a subgroup $\Gamma\subset\hbox{Aut}(X)$ {\it hereditarily ergodic}
if for every closed connected virtually $\Gamma$-invariant subgroup $Y$ of $X$, the
 action of $\Gamma_Y$ on $Y$ is ergodic.

It is not hard to check that for abelian groups of automorphisms of compact
connected finite-dimensional group $X$, the notions of ergodicity and
hereditary ergodicty coincide (this fails, in general, for infinite-dimensional
groups $X$
--- see Example \ref{ex:her} below). Hence, Berend's theorem in this case
states that hereditary ergodicity is equivalent to existence of
an automorphism which is
mixing of all orders.
The following theorem generalizes this result to solvable
groups of automorphisms.

\begin{thm}\label{th:solv1}
Let $X$ be a compact connected finite-dimensional abelian group and
$\Gamma$ a solvable subgroup of automorphisms of $X$. Then the following
statements are equivalent:
\begin{enumerate}
\item[(a)] $\Gamma$ is hereditarily ergodic.
\item[(b)]  $\Gamma$ contains an abelian subgroup which is mixing of all orders.
\end{enumerate}
\end{thm}

Note that the assumption in Theorem \ref{th:solv1} that the group
$\Gamma$ is solvable is essential (see Example \ref{ex:not_erg} below).
Also, Theorem \ref{th:solv1} fails without the assumption that $X$ is
finite-dimensional (see Examples \ref{ex:her} and \ref{ex:last} below).

According to the {\it Rosenblatt's alternative} (see
\cite{ros}),  any finitely
generated solvable group is
either virtually nilpotent or contains a free nonabelian subsemigroup.
In the latter case Theorem \ref{th:solv1} can be strengthened as follows:

\begin{thm}\label{th:solv2}
Let $X$ be a compact connected finite-dimensional abelian group and
$\Gamma$ a solvable group of automorphisms of $X$, which is not
virtually nilpotent. Then the following
statements are equivalent:
\begin{enumerate}
\item[(a)] $\Gamma$ is hereditarily ergodic.
\item[(b)]  $\Gamma$ contains a free nonabelian subsemigroup which is mixing of all orders.
\end{enumerate}
\end{thm}

Combining Theorems \ref{th:solv1} and \ref{th:solv2}, we deduce

\begin{cor}\label{c:irr}
Let $X$ be a compact connected finite-dimensional abelian group and
$\Gamma$ a solvable strongly irreducible group of automorphisms of
$X$. Then $\Gamma$ contains an abelian subgroup which is mixing of all orders.
Moreover, if $\Gamma$ is not virtually nilpotent, then $\Gamma$
contains a free nonabelian subsemigroup which is mixing of all orders.
\end{cor}

Without the assumption that the group $\Gamma$ is solvable, Corollary
\ref{c:irr} fails (see Example \ref{ex:not_erg} below).

It follows from the {\it Tits alternative} (see \cite{tits} or
\cite[Section 5J]{mor}) that any finitely generated 
subgroup of $\hbox{Aut}(X)$ is either virtually solvable or contains a
nonabelian free group. Recently, E.~Breuillard and  T.~Gelander
proved a {\it topological Tits alternative} (see \cite{br_g}):
any finitely generated matrix group either contains a Zariski open solvable
subgroup or a Zariski dense free subgroup.
Utilizing this result, we obtain

\begin{thm}\label{th:tits}
Let $X$ be a compact connected finite dimensional abelian group and
$\Gamma$ an ergodic (hereditarily ergodic, strongly irreducible)
subgroup of $\hbox{\rm Aut}(X)$ which is not virtually
solvable. Then $\Gamma$ contains a free
nonabelian ergodic (hereditarily ergodic, strongly irreducible) subgroup.
\end{thm}

Example \ref{ex:not_erg} below illustrates that an ergodic group may contain no ergodic elements.

\subsection{}
Some special cases of the above results appeared in \cite{bg}.
Note that in \cite{bg} we used a slightly different definition for mixing
(borrowed from \cite{ber}), but in this paper we adopt the definition
from \cite{sch}. The relation between these two definitions is
quite straightforward (see Corollary \ref{p:strong}).

The paper is organized as follows.
The main theorem (Theorem \ref{th:general}) is proved in Section
\ref{sec:proof}. The rest of the results stated in Subsection
\ref{sec:spectrum} are proved in Section \ref{sec:cor}.
The results about mixing groups of automorphisms (stated in Subsection
\ref{subsec:group}) are proved in Section \ref{sec:group}.
The theorems from Subsection \ref{s:e} are proved in Section \ref{sec:erg}.
Section \ref{sec:ex} contains some examples and counterexamples related
to the results of this paper.

\subsection{Acknowledgement}
We would like to thank Y.~Benoist and H.~Oh for helpful discussions
and D.~Berend and T.~Ward for useful comments about this paper.

\section{Mixing and linear relations 
(proof of Theorem \ref{th:general})}\label{sec:proof}

Let $X$ be a compact connected abelian group of finite dimension $d$. We identify its character group $\hat X$
with a subgroup of $\mathbb{Q}^d$. Then every endomorphism $T$ of $X$ induces a linear map $\hat T$ of $\mathbb{Q}^d$.

We recall the well-known characterization of mixing:

\begin{lem}\label{l:mix_cond}
The set $\{T_1,\ldots, T_s\}$ of epimorphisms of $X$ is  mixing iff
there are no $x_1,\ldots,x_{s}\in \mathbb{Q}^d$ such that $(x_1,\ldots, x_{s})\ne (0,\ldots, 0)$
and for infinitely many $n\ge 1$,
$$
\hat T_1^nx_1+\cdots+ \hat T_s^n x_s=0.
$$
\end{lem}


As an application of Lemma \ref{l:mix_cond}, we show that
the set of epimorphisms $\{T_1,\ldots,T_{s}\}$, $s>d$, is not mixing
provided that the linear
maps $\hat T_1,\ldots,\hat T_{s}$ have the same characteristic
polynomial. 

\begin{pro}\label{p_not_mix}
Let $T_1,\ldots, T_s$ be epimorphisms of $X$ 
and assume that there exists a polynomial $p(x)\in\mathbb{Q}[x]$ with
$\deg p<s$ such that $p(\hat T_i)=0$ for $i=1,\ldots, s$.
Then $\{T_1,\ldots, T_s\}$ is not mixing. 
\end{pro}

\begin{proof}
To prove the proposition, it suffices to construct $(x_1,\ldots,x_{s})\in (\mathbb{Q}^d)^{s}-\{(0,\ldots,0)\}$ such that
\begin{equation}\label{eq:not_j_m3}
\hat T_1^nx_1+\cdots +\hat T_{s}^n x_{s}=0
\end{equation} 
for infinitely many $n$.

Let $p(x)=p_1(x)^{e_1}\cdots p_l(x)^{e_l}$ where $p_i(x)\in\mathbb{Q}[x]$,
$i=1,\ldots,l$, are distinct and irreducible. Let $d_i=\deg p_i$ and $\lambda_{i,j}$, $j=1,\ldots,d_i$,
be the roots of $p_i$. Let $P_{i,j,k}\in \hbox{M}(d,\mathbb{Q}(\lambda_{i,j}))$ be the projection on
the root space of $T_k$ corresponding to $\lambda_{i,j}$. Then
$$
\hat T_k^nP_{i,j,k}=\lambda_{i,j}^n\sum_{u=0}^{e_i-1} n^u A_{i,j,k,u}
$$
for some $A_{i,j,k,u}\in \hbox{M}(d,\mathbb{Q}(\lambda_{i,j}))$. Since the coefficients of the $(e_id)\times sd$
matrix
$$
B_{i,j}=(A_{i,j,k,u}: u=0,\ldots,e_i-1,\, k=1,\ldots,s)
$$
lie in $\mathbb{Q}(\lambda_{i,j})$, we have
$$
\hbox{rank}_\mathbb{Q} (B_{i,j})\le e_i d\cdot [\mathbb{Q}(\lambda_{i,j}):\mathbb{Q}]=e_idd_i.
$$
It follows that for the $(le_id)\times sd$ matrix
$$
C=(B_{i,1}:i=1,\ldots l),
$$
$$
\hbox{rank}_\mathbb{Q} (C)\le \sum_{i=1}^l e_idd_i=(\deg p)\cdot d.
$$
Hence, there exists a vector
$x=(x_1,\ldots,x_s)\in (\mathbb{Q}^d)^s-\{(0,\ldots,0)\}$ such that $Cx=0$. Then
\begin{equation}\label{eq:last}
\sum_{k=1}^s \hat T_k^n P_{i,1,k} x_k=0
\end{equation}
for every $i=1,\ldots,l$. For fixed $i,k$, the Galois group $\hbox{Gal}(\mathbb{C}/\mathbb{Q})$
permutes transitively the roots $\lambda_{i,j}$ and the matrices $P_{i,j,k}$, $j=1,\ldots,e_i$.
Hence, if follows from (\ref{eq:last}) that for every $j=1,\ldots,e_i$,
\begin{equation}\label{eq:last2}
\sum_{k=1}^s \hat T_k^n P_{i,j,k} x_k=0.
\end{equation}
Summing (\ref{eq:last2}) over $i$ and $j$, we deduce (\ref{eq:not_j_m3}).
\end{proof}

To analyze the equation in Lemma \ref{l:mix_cond}, we use the
following statement sometimes referred to as Kronecker's lemma (see
\cite[p. 27]{ew}):

\begin{lem}[Kronecker]\label{l:abs_1}
If $\lambda$ is an algebraic integer such that all of its conjugates have absolute value one,
then $\lambda$ is a root of unity.
\end{lem}

We also mention an equivalent formulation of Kronecker's lemma, which
we use latter: if $\lambda$ is an element of a number field $K$ and $\lambda$ is
not a root of unity, then there exists an absolute value $|\cdot|_v$ on $K$ such
that $|\lambda|_v\ne 1$.

Note that $\{T_1,\ldots,T_s\}$ is mixing iff $\{T_1^l,\ldots,T_s^l\}$ is 
mixing for some (all) $l\ge 1$. This observation implies that in the
proof of Theorem \ref{th:general} we may assume
without loss of generality that 
\begin{equation}\label{eq:spec2}
\lambda,\mu\in \bigcup_{k=1}^{s}\hbox{Spec}(T_k)\hbox{ and
}\lambda^{-1}\mu\hbox{ a root of unity }\, \Rightarrow\, \lambda=\mu. 
\end{equation}
Under this assumption, Theorem \ref{th:general} can be restated as follows

\begin{thm}\label{th:general2}
Let $T_1,\ldots,T_s$ be epimorphisms of $X$ that satisfy
(\ref{eq:spec2}). Then the set $\{T_1,\ldots,T_s\}$ is mixing iff
for every subset $\{k_1,\ldots,k_r\}\subset\{1,\ldots,s\}$ and every
$\lambda\in\mathbb{C}$,
there are no $\lambda$-eigenvectors of $\hat T_{k_1},\ldots, \hat T_{k_r}$
that are linearly dependent over $\mathbb{Q}$.
\end{thm}  

\begin{proof}
Suppose that there exist a nonempty subset $S\subset \{1,\ldots,s\}$,
$\alpha_k\in\mathbb{Q}-\{0\}$, $k\in S$, and eigenvalues $w_k$ for
$\hat T_k$, $k\in S$, with the same eigenvalue $\lambda$ such that
$$
\sum_{k\in S} \alpha_k w_k=0.
$$
This implies that the subspace
$$
V=\left\{(v_k)\in (\mathbb{C}^d)^{|S|}:\; \sum_{k\in S} \alpha_k
  \hat T_k^n v_k=0\;\hbox{ for all $n\ge 1$}\right\}
$$
is not trivial. Since this subspace is defined over $\mathbb{Q}$, it
contains a nonzero rational vector $(x_k:k\in S)$
that gives a nonzero solution of the equation
\begin{equation}\label{eq:mix3}
\sum_{k=1}^{s} \hat T_k^nx_k=0.
\end{equation}
Hence, by Lemma \ref{l:mix_cond}, the set $\{T_1,\ldots, T_s\}$ is not mixing.

Conversely, suppose that the set $\{T_1,\ldots, T_s\}$ is not mixing.
Then by Lemma \ref{l:mix_cond}, there exists $(x_1,\ldots,x_{s})\in
(\mathbb{Q}^d)^{s}-\{(0,\ldots,0)\}$ such that (\ref{eq:mix3})
holds for infinitely many $n\ge 1$.

Let 
$$p_k(x)=p_{k,1}(x)^{m_{k,1}}\cdots p_{k,l_k}(x)^{m_{k,l_k}}$$
be the characteristic polynomial of $\hat T_k$, $k=1,\ldots,s$, where $p_{k,i}(x)\in\mathbb{Q}[x]$ are distinct and irreducible over $\mathbb{Q}$. Let $d_{k,i}=\deg(p_{k,i})$.
For a root $\lambda$ of $p_k$, denote by $V_{k}^\lambda$
the root subspace of $T_k$ with respect to $\lambda$. Then
$$
\mathbb{C}^d=\bigoplus_{\lambda:\, p_k(\lambda)=0} V_{k}^\lambda.
$$
Note that for fixed $k$ and $i$, the Galois group $\hbox{Gal}(\mathbb{C}/\mathbb{Q})$ permutes transitively the spaces $V_{k}^\lambda$
where $\lambda$ satisfies $p_{k,i}(\lambda)=0$. This implies that the subspaces
$$
V_{k,i}=\bigoplus_{\lambda:\,p_{k,i}(\lambda)=0} V_{k}^\lambda
$$
are rational. Then
\begin{equation}\label{eq:vki}
\mathbb{Q}^d=\bigoplus_{i=1}^{l_k} V_{k,i}(\mathbb{Q}),
\end{equation}
and there exist vectors $x_{k,i}\in V_{k,i}(\mathbb{Q})$, not all zero,
such that
\begin{equation}\label{eq:mix4}
\sum_{k=1}^{s}\sum_{i=1}^{l_k} \hat T_k^nx_{k,i}=0
\end{equation}
for infinitely many $n\ge 1$.
For a root $\lambda$ of $p_{k,i}$,
let $P_{k}^\lambda$ denote the projection from $V_{k,i}$ on the root space $V_{k}^\lambda$.
Since
$$
\hat T_k|_{V_{k}^\lambda}=\lambda(id+N_{k}^\lambda)
$$
where $N_{k}^\lambda:V_{k}^\lambda\to V_{k}^\lambda$ is nilpotent linear map such that $(N_{k}^\lambda)^{m_{k,i}}=0$,
we have
$$
\hat T_k^nP_k^\lambda=\lambda^n\sum_{u=0}^{m_{k,i}-1} {n\choose u} A_{k,u}^\lambda
$$
where $A_{k,u}^\lambda:V_{k,i}\to V_{k}^\lambda$ are linear maps and $A_{k,0}^\lambda=P_{k}^\lambda$.
With respect to a rational basis on
$V_{k,i}$, $A_{k,u}^\lambda$ is represented by $d\times (\dim V_{k,i})$ matrix with coefficients in $\mathbb{Q}(\lambda)$. Then (\ref{eq:mix4}) is equivalent to
\begin{equation}\label{eq:mix4_1}
\sum_{k=1}^{s}\sum_{i=1}^{l_k}\sum_{\lambda:\, p_{k,i}(\lambda)=0}
\sum_{u=0}^{m_{k,i}-1} \lambda^n {n\choose u} A_{k,u}^\lambda x_{k,i}=0.
\end{equation}

Denote by $K$ the number field generated by the eigenvalues of $T_i$,
$i=1,\ldots,s$, and let $\mathcal{V}_K$ be the set of absolute values of $K$.

Since (\ref{eq:mix4_1}) holds for infinitely many $n$, it is equivalent to the system of equations
\begin{equation}\label{eq:mix5}
\sum_{k,i,\lambda,j}^{} {}^\prime \sum_{u=0}^{m_{k,i}-1}\lambda^n
{n\choose u} A_{k,u}^\lambda x_{k,i}=0, \quad\delta>0,
\end{equation}
where the sum $\sum'$ is taken over those $\lambda$'s such that
$p_{k,i}(\lambda)=0$ and $|\lambda|_v=\delta$,
$v\in\mathcal{V}_K$. 
Conjugating (\ref{eq:mix5}) by $\sigma\in\hbox{Gal}(\mathbb{C}/\mathbb{Q})$, we deduce
that (\ref{eq:mix5}) is equivalent to the system of equations
\begin{equation}\label{eq:mix6}
\sum_{k,i,\lambda,j}{}{}^{''} \sum_{u=0}^{m_{k,i}-1} \lambda^n {n\choose u}
A_{k,u}^\lambda x_{k,i}=0,\quad\delta_{\sigma,v}>0,\, \sigma\in
\hbox{Gal}(\mathbb{C}/\mathbb{Q}),\, v\in\mathcal{V}_K,
\end{equation}
where the sum $\sum''$ is taken over $\lambda$'s such that $p_{k,i}(\lambda)=0$ and $|\lambda^\sigma|_v=\delta_{\sigma,v}$ for
every $\sigma\in \hbox{Gal}(\mathbb{C}/\mathbb{Q})$ and $v\in\mathcal{V}_K$.

If 
$$
\lambda,\mu\in\bigcup_{k=1}^{s} \hbox{Spec}(T_k)
$$
and $|\lambda^\sigma|_v=|\mu^\sigma|_v$ for every
$\sigma\in \hbox{Gal}(\mathbb{C}/\mathbb{Q})$ and $v\in\mathcal{V}_K$,
then $\lambda^{-1}\mu$ is a root of unity by Lemma \ref{l:abs_1}, and by (\ref{eq:spec2}),
$\lambda=\mu$. Hence, (\ref{eq:mix6}) is equivalent to the system of equations
\begin{equation}\label{eq:mix7}
\sum_{k,i:\, p_{k,i}(\lambda)=0}\sum_{u=0}^{m_{k,i}-1}{n\choose u} A_{k,u}^\lambda x_{k,i}=0, \quad\lambda\in\bigcup_{k=1}^{s} \hbox{Spec}(T_k).
\end{equation}
Let $m_\lambda=\max\{m_{k,i}:p_{k,i}(\lambda)=0\}$. Since (\ref{eq:mix7}) holds for infinitely many $n$,
it is equivalent to
\begin{equation}\label{eq:mix8}
\sum_{k,i:\, p_{k,i}(\lambda)=0} A_{k,u}^\lambda x_{k,i}=0, \quad\lambda\in\bigcup_{k=1}^{s} \hbox{Spec}(T_k),\, u=0,\ldots,m_\lambda-1.
\end{equation}
For every $k=1,\ldots,s$ and $i=1,\ldots, l_k$, choose $\lambda_{k,i}$ such that $p_{k,i}(\lambda_{k,i})=0$.
If $p_{k,i}$'s have a common root for different $k$'s, we choose the same $\lambda_{k,i}$.
Let 
$$
\Lambda=\{\lambda_{k,i}:k=1,\ldots, s,\, i=1,\ldots, l_k\}.
$$
Note that for $\sigma\in\hbox{Gal}(\mathbb{C}/\mathbb{Q})$, we have 
$$
\sigma(V_k^\lambda)=V_k^{\sigma(\lambda)},\;\sigma(P_{k}^\lambda)=P_k^{\sigma(\lambda)},\; \sigma(N_k^\lambda)=N_k^{\sigma(\lambda)},\;\sigma(A_{k,u}^\lambda)=A_{k,u}^{\sigma(\lambda)}.
$$
Since the polynomial $p_{k,i}$ is irreducible,
the Galois group $\hbox{Gal}(\mathbb{C}/\mathbb{Q})$ acts transitively on the set of roots of $p_{k,i}$.
Hence, if  (\ref{eq:mix8}) holds for $\lambda=\lambda_{k,i}$, then it holds for all $\lambda$'s such that
$p_{k,i}(\lambda)=0$. Therefore, (\ref{eq:mix8}) is equivalent to
\begin{equation}\label{eq:mix9}
\sum_{k,i:\, p_{k,i}(\lambda)=0}A_{k,u}^\lambda x_{k,i}=0, \quad\lambda\in \Lambda,\, u=0,\ldots,m_\lambda-1.
\end{equation}
Since polynomials $p_{k,i}$, $i=1,\ldots, l_k$, have no common roots, it follows that for every $k=1,\ldots, s$
 and $\lambda\in\Lambda$, there is at most one $i$ such that $p_{k,i}(\lambda)=0$.
Hence, the system of equations (\ref{eq:mix9}) splits into independent systems of equations
\begin{equation}\label{eq:mix10}
\sum_{k,i:\, p_{k,i}(\lambda)=0}A_{k,u}^\lambda x_{k,i}=0, \quad u=0,\ldots,m_\lambda-1
\end{equation}
indexed by $\lambda\in \Lambda$.  Therefore, (\ref{eq:mix9}) has a nontrivial solution iff for some $\lambda\in \Lambda$,
(\ref{eq:mix10}) has a nontrivial solution.

Let $\lambda\in \Lambda$ be such that (\ref{eq:mix10}) has a
nontrivial solution and $u_0\in\{0,\ldots,m_\lambda-1\}$ be maximal index such
that (\ref{eq:mix10}) contains nonzero terms.
Since
$$
A^\lambda_{k,u}=\lambda^{-u}(\hat T_k-\lambda)^u|_{V_k^\lambda},\quad u\ge 1,
$$
it follows that (nonzero) vectors $A_{k,u_0}^\lambda x_{k,i}$ are
eigenvectors of $\hat T_k$ with eigenvalue $\lambda$ which are
linearly dependent over $\mathbb{Q}$. This proves the theorem.
\end{proof}

\section{Proofs of Corollaries formulated in Subsection \ref{sec:spectrum}}\label{sec:cor}

\begin{proof}[Proof of Corollary \ref{cor:com}]
It is clear that (a)$\Rightarrow$(b).

Suppose that $\hat T_i^{-1}\hat T_j$ has a root of unity as an
eigenvalue. Since $\hat T_i$ and $\hat T_j$ commute, this implies that for some $l\ge 1$, the subspace
$$
V=\{v\in \mathbb{C}^d:\, \hat T_i^l v= \hat T_j^l v\}
$$
is not $\{0\}$. Since $V$ is $\{\hat T_i^l,\hat T_j^l\}$-invariant and
$\hat T_i^l|_V=\hat T_j^l|_V$, it follows that $\hat T_i^l$ and $T_j^l$ have common eigenvector
with the same eigenvalue. Hence, by Theorem \ref{th:general}, $\{T_i, T_j\}$
is not mixing. This shows that (b)$\Rightarrow$(c).

To prove that (c)$\Rightarrow$(a), suppose that (c) holds, but $\{T_1,\ldots, T_s\}$ is not mixing.
Then by Theorem \ref{th:general}, there exist a nonempty $S\subset \{1,\ldots,s\}$,
$\alpha_k\in\mathbb{Q}-\{0\}$, $k\in S$,
and eigenvectors $w_k\in \mathbb{C}^d$ of $\hat T_k^l$, $k\in S$, with the same
eigenvalue $\lambda$ such that
$$
\sum_{k\in S} \alpha_k w_k=0.
$$
Hence, we have a nonzero vector space
$$
V=\left\{(v_k)\in (\mathbb{C}^d)^{|S|}:\, \sum_{k\in S} \alpha_k v_k=0,\;\; \hat
T_k^l v_k=\lambda v_k\;\hbox{ for $k\in S$}\right\}.
$$
Since $\hat T_i$'s commute, this vector space is $\{\hat T_1^l,\ldots,
\hat T_s^l\}$-invariant, and it contains a common eigenvector
$v=(v_k:k\in S)$:
$$
\hat T_k^lv=\lambda_k v,\quad k\in S.
$$
Let $S_0=\{k\in S:\, v_k\ne 0\}$. Note that $|S_0|>1$. For $k\in S_0$,
$\lambda_k=\lambda$. Hence, $\hat T_i^{-l}\hat T_j^l=(\hat
T_i^{-1}\hat T_j)^l$ has eigenvalue 1  for $i,j\in S_0$. This
contradicts (c). Hence, (c)$\Rightarrow$(a).
\end{proof}

\begin{proof}[Proof of Corollary \ref{th:two}]
If $\{T_1,T_2\}$ is mixing on $X$, then clearly, $\{T_1^l,T_2^l\}$ is mixing on $X/Y$.
Hence, one direction of the corollary is obvious.

Suppose that $\{T_1,T_2\}$ is not mixing. Then, by Theorem \ref{th:general},
$\hat T_1^l$ and $\hat T_2^l$ have common eigenvector with the same
eigenvalue for some $l\ge 1$.
This eigenvector is contained in the rational subspace 
$$
V=\{v\in \mathbb{C}^d:\, \hat T_1^lv=\hat T_2^lv,\, \hat T_1^l\hat T_2^lv=\hat T_2^l\hat T_1^lv\}.
$$
Consider the subgroup
$$
Y=\{x\in X:\, \chi(x)=1\quad\hbox{for}\;\; \chi\in V\cap \hat X\}.
$$
Since $V$ is rational, $V\cap \hat X\ne 0$ and $Y\ne X$, and since the subspace $V$ is $\{\hat T_1^l,\hat T_2^l\}$-invariant,
the subgroup $Y$ is $\{T_1^l,T_2^l\}$-invariant. 
The character group of $X/Y$ is $V\cap \hat X$,
and $\hat T_1^l=\hat T_2^l$ on $V\cap \hat X$. Hence, it follows that $T_1^l=T_2^l$ on $X/Y$.
This proves the corollary.
\end{proof}

\begin{proof}[Proof of Corollary \ref{th:cond_gen}]
Under the assumption in (a), for every $\lambda\in \mathbb{C}$, there
is at most one $\hat T_i^l$ with the eigenvalue $\lambda$. Hence,
the maps $\hat T_i^l$ cannot have linearly dependent eigenvectors with the
same eigenvalue, and by Theorem \ref{th:general}, the set
$\{T_1,\ldots, T_s\}$ is mixing.

Suppose that there exist $S\subset \{1,\ldots,s\}$ with $|S|=r>d$, $l\ge 1$, and
$\lambda\in\mathbb{C}$ such that
$$
\lambda\in \hbox{Spec}(\hat T_k^l)\quad\hbox{ for $k\in S$.}
$$
We are going to show now that the set $\{T_k^l:\, k\in S\}$ is not mixing.
This will imply that the set $\{T_1,\ldots, T_s\}$ is not mixing as well.

Denote by $q(x)\in \mathbb{Q}[x]$ the minimal polynomial of $\lambda$ and
consider a rational subspace
$$
W_{k}=\{v\in \mathbb{C}^d:\, q(\hat T_k^l)v=0\}.
$$
Note that $W_k$ contains all $\mu$-eigenspaces of $\hat T_k^l$ such that $q(\mu)=0$.
In particular,
$$
\dim W_k\ge \deg (q).
$$
Denote by $P_k$ the projection from $W_k$ to the $\lambda$-eigenspace
of $\hat T_k^l$. 
According to Theorem \ref{th:general}, it suffices to show that
there exist $x_{k}\in W_k(\mathbb{Q})$, not all zero, such that
\begin{equation}\label{eq:mix11}
\sum_{k\in S}P_k x_{k}=0.
\end{equation}
Choose a rational basis in $W_{k}$.
With respect to this basis,  the linear map $P_k$
is represented by $d\times (\dim W_k)$ matrix with coefficients in $\mathbb{Q}(\lambda)$.
Consider the  $d\times (\dim W^\lambda_k\cdot |S|)$-matrix 
$$
P=(P_{k}: k\in S).
$$ 
Since the coefficients of $P$ are in $\mathbb{Q}(\lambda)$, 
$$
\hbox{rank}_\mathbb{Q} (P)\le d\cdot[\mathbb{Q}(\lambda):\mathbb{Q}]=d\cdot \deg (q)<|S|\cdot\dim W_k.
$$
Hence, there exists nonzero vector 
$$
x=(x_{k}:k\in S)\in \prod_{k\in S} W_k(\mathbb{Q})\simeq\mathbb{Q}^{|S|\cdot\dim W_k}
$$
such that $P\cdot x=0$. Hence, (\ref{eq:mix11}) has a nonzero
solution. This proves (b).

Now we prove (c). Suppose that for every $l\ge 1$ and $S\subset \{1,\ldots,s\}$ such that  $|S|>d$,
$$
\bigcap_{i\in S} \hbox{\rm Spec}(T_i^l) = \emptyset,
$$
Then if $\hat T_i^l$'s have linearly dependent
eigenvectors with the same eigenvalue,  there is a subset of $\hat
T_i^l$'s of cardinality at most $d$ with the same property. Hence, it
follows from Theorem \ref{th:general} that $\{T_1,\ldots, T_s\}$
is mixing iff every subset of cardinality $d$ is mixing.
\end{proof}

\begin{proof}[Proof of Corollary \ref{thm:limit}]
We choose $l\ge 1$ so that $T_1^l,\ldots, T_s^l$ satisfy condition
(\ref{eq:spec2}). It suffices to prove the corollary for $k=0$, and
to simplify calculations, we also assume that $l=1$. The proof of the
general case easily reduces to this situation.

We use the notation introduced in the proof of Theorem \ref{th:general2}.

For $\chi_1,\ldots,\chi_s\in \hat X$,
\begin{equation}\label{eq:char}
\int_X \chi_1(T_1^{n}x)\cdots \chi_s(T_s^{n}x)\,dm(x)=\left\{
\begin{tabular}{l} $1$ if $\hat T_1^{n}\chi_1+\cdots+\hat
  T_s^{n}\chi_s=0$,\\
$0$ if $\hat T_1^{n}\chi_1+\cdots+\hat
  T_s^{n}\chi_s\ne 0$.
\end{tabular}
\right.
\end{equation}
Denote by $Q_{k,i}$ the projection on the space $V_{k,i}$ with respect to
the decomposition (\ref{eq:vki}).
If
\begin{equation}\label{eq:eq0}
\hat T_1^{n}\chi_1+\cdots+\hat
  T_s^{n}\chi_s=0
\end{equation}
for infinitely many $n$, then by the proof of Theorem \ref{th:general},
\begin{equation}\label{eq:eq1}
\sum_{k,i: p_{k,i}(\lambda)=0} A_{k,u}^\lambda Q_{k,i}
\chi_k=0\quad\hbox{for $\lambda\in\Lambda$, $u=0,\ldots, m_\lambda-1$.}
\end{equation}
Conversely, (\ref{eq:eq1}) implies that (\ref{eq:eq0}) holds for every
$n\ge 1$. Denote by $\Delta$ the set of $(\chi_1,\ldots, \chi_s)\in
\hat X^s$ such that (\ref{eq:eq1}) holds. We claim that
for every $f_1,\ldots, f_s\in L^\infty(X)$,
$$
\lim_{n\to\infty} \int_X f_1(T_1^{n}x)\cdots
f_s(T_s^{n}x)\,dm(x)=\sum_{(\chi_1,\ldots,\chi_s)\in\Delta}
\hat f_1(\chi_1)\cdots \hat f_s(\chi_s).
$$
When $f_1,\ldots, f_s$ are
characters, this follows from (\ref{eq:char}).
For general $L^\infty$-functions, the claim is proved by the standard
approximation argument. 
\end{proof}

\begin{proof}[Proof of Corollary \ref{th:ex}(a)]
Let us choose linearly dependent over $\mathbb{Q}$ vectors $v_1,\ldots,v_s\in
 \mathbb{Z}^d$ such that every proper subset of $\{v_1,\ldots,v_s\}$ is linearly independent over
$\mathbb{Q}$.
There exist nilpotent matrices $N_1,\ldots,
N_s\in\hbox{M}(d,\mathbb{Z})$ such that
$$
\hbox{Ker}(N_i)=\left<v_i\right>.
$$
Set $T_i=id+N_i$. Then it follows from Theorem \ref{th:general} that
the set $\{T_1,\ldots,T_s\}$ is not mixing, but its every proper
subset is mixing.
\end{proof}

To prove Corollary \ref{th:ex}(b), we need a lemma:

\begin{lem}\label{lem:polynom}
For every $d\ge 1$, there exists an irreducible monic polynomial $p(x)\in\mathbb{Z}[x]$ which
has real roots with different absolute values.
\end{lem}

\begin{proof}
Consider the polynomial
$$
p(x)=(x-q)\cdots (x-dq)+q
$$
where $q$ is a prime number. Note that this polynomial is irreducible
by the Eisenstein criterion (see \cite[IV \S 3]{lang}). Let us assume that $d$ is even (the
argument for odd $d$ is analogous). Then for sufficiently large $q$, we have
\begin{align*}
p((4i+1)q/2)&\ge (q/2)^d+q>0,\quad i=0,\ldots, d/2,\\
p((4i+3)q/2)&\le -(q/2)^d+q<0,\quad i=0,\ldots, d/2-1.
\end{align*}
This implies that $p(x)$ has $d$ distinct positive real roots.
\end{proof}

\begin{proof}[Proof of Corollary \ref{th:ex}(b)]
Let $p(x)$ be as in Lemma \ref{lem:polynom} and let $T_1\in\hbox{M}(d,\mathbb{Z})$ 
has $p(x)$ as its characteristic polynomial.
Denote by $\lambda_i$, $i=1,\ldots,d$, the roots of $p(x)$ and
by $\sigma_i$ the embedding $\mathbb{Q}(\lambda_1)\to\mathbb{R}$
such that $\lambda_1\mapsto\lambda_i$.
Let $\{v_1,\ldots,v_d\}$ be an integral basis of $\mathbb{Q}(\lambda_1)$. It is well-known
that
\begin{equation}\label{eq:det}
\det (v_{k}^{\sigma_i}: i,k=1,\ldots,d)\ne 0.
\end{equation}
Let 
$$
A=\hbox{diag}(1,\ldots,1,2,\ldots,s-1).
$$
Note that $A$ has minimal polynomial 
$$q(x)=\prod_{l=1}^{s-1}(x-l)= \sum_{j=1}^{s} q_jx^{j-1}$$
and $q_j\ne 0$ for all $j$.
Put
$$
w_j=A^{j-1}\cdot{}^t(v_{1},\ldots,v_{d}),\quad j=1,\ldots, d.
$$
It follows from (\ref{eq:det}) that $w_j^{\sigma_i}$, $i=1,\ldots, d$,
are linearly independent over $\mathbb{Q}$.
The Galois group
$\hbox{Gal}(\mathbb{C}/\mathbb{Q})$ permutes the vectors $w_j^{\sigma_i}$, $i=1,\ldots, d$.
Therefore,
\begin{equation}\label{eq:Q}
\mathbb{Q}^d=\left\{\sum_{i=1}^d a^{\sigma_i}w_j^{\sigma_i}: a\in\mathbb{Q}(\lambda_1)\right\}.
\end{equation}
for every $j=1,\ldots,d$. Define $T_j\in\hbox{M}(d,\mathbb{R})$ such that
$$
T_jw_j^{\sigma_i}=\lambda_i w_j^{\sigma_i},\quad i=1,\ldots,d.
$$
Then $\det (T_j)\ne 0$ and $T_j^\sigma=T_j$ for every $\sigma\in\hbox{Gal}(\mathbb{C}/\mathbb{Q})$.
Thus, $T_j\in\hbox{GL}(d,\mathbb{Q})$. Multiplying $T_j$'s and $\lambda_i$'s by an integer we may assume
that $T_j$'s have integer entries.


We claim that there is $(x_1,\ldots,x_{s})\in (\mathbb{Q}^d)^{s}-\{(0,\ldots,0)\}$
such that
\begin{equation}\label{eq:mix1}
\sum_{j=1}^{s} T^n_jx_j=0
\end{equation}
for every $n\ge 1$,
and for any $J\subsetneq \{1,\ldots, s\}$, there is no $(x_j:j\in J)\in (\mathbb{Q}^d)^{|J|}-\{(0,\ldots,0)\}$
such that
\begin{equation}\label{eq:mix2}
\sum_{j\in J}T^n_jx_j=0
\end{equation}
for infinitely many $n$.
By Lemma \ref{l:mix_cond}, this implies that $\{T_1,\ldots,T_s\}$ is
not mixing, but its every proper subset is mixing. 

Put $x_j=q_j\sum_{i=1}^d w_j^{\sigma_i}\in\mathbb{Q}^d$. Note that
$$
\sum_{i=1}^d w_j^{\sigma_i}=A^j\left(\sum_{i=1}^d w_1^{\sigma_i}\right)\ne 0
$$
by (\ref{eq:det}). We have
$$
\sum_{j=1}^{s} T^n_jx_j=\sum_{i=1}^d \sum_{j=1}^s q_j\lambda_i^{n}w_j^{\sigma_i}
=\sum_{i=1}^d \lambda_i^{n}\sum_{j=1}^s q_jA^{j-1}w_1^{\sigma_i}=0.
$$
This proves (\ref{eq:mix1}).

It follows from  (\ref{eq:Q}) that  (\ref{eq:mix2}) is equivalent to existence of 
$$
(a_j:j\in J)\in \mathbb{Q}(\lambda_1)^{|J|}-\{(0,\ldots,0)\}
$$
such that
$$
\sum_{j\in J}\sum_{i=1}^d T^n_ja_j^{\sigma_i}w^{\sigma_i}_j=0
$$
for infinitely many $n$. Then
$$
\sum_{i=1}^d \lambda^n_i\sum_{j\in J}a_j^{\sigma_i}w^{\sigma_i}_j=0.
$$
Since $\lambda_i$'s have different absolute values, this implies that
\begin{equation}\label{eq:rat}
\sum_{j\in J} a_jw_j=\left(\sum_{j\in J} a_jA^{j-1}\right)w_1=0.
\end{equation}
Let 
$$
r(x)=\sum_{j\in J} a_j x^{j-1}=b\prod_{j=1}^l (x-\mu_j).
$$
for $l\le s-1$ and $b,\mu_j\in \mathbb{C}$. We have
$r(A)w_1=0$ and $\tilde{r}(A)w_1=0$ with
$$
\tilde{r}(x)=\prod (x-\mu_j)
$$
where the product is taken over $\mu_j$ which are eigenvalues of $A$. In particular, $\tilde{r}(x)\in\mathbb{Q}[x]$.
Then $\tilde{r}(A)w_1^{\sigma_i}=0$ for every $i=1,\ldots, d$. It follows from (\ref{eq:det}) that
$\tilde{r}(A)=0$. Since the minimal polynomial $q(x)$ of $A$ has degree $s-1$, this implies that $r(x)$ is a scalar multiple of
$q(x)$. In particular, $a_j\ne 0$ for every $j=1,\ldots, s$, which is a contradiction.
\end{proof}

\begin{proof}[Proof of Corollary \ref{p:strong}]
It is clear that if the set $\{T_1,\ldots,T_{s-1},id\}$ is mixing,
then every $T_i$ is mixing and $\{T_1,\ldots, T_{s-1}\}$ is mixing as
well.

Conversely, suppose that the set $\{T_1,\ldots,T_{s-1},id\}$ is not
mixing. Then for some $l\ge 1$, the linear maps $\hat
T_1^l,\ldots,\hat T_{s-1}^l,id$ have linearly dependent eigenvectors
with the same eigenvalue $\lambda$. Note that if $\lambda=
1$, then some $\hat T_k^l$ has eigenvalue one, and $T_k$ is not mixing.
Otherwise, it follows that the linear maps $\hat
T_1^l,\ldots,\hat T_{s-1}^l$ have linearly dependent eigenvectors
with the same eigenvalue. Hence, by Theorem \ref{th:general}, the 
set $\{T_1,\ldots, T_{s-1}\}$ is not mixing.
\end{proof}

\section{Mixing groups and semigroups}\label{sec:group}

\begin{pro}\label{p:mix}
Let $X$ be any compact abelian group and $\Gamma$
a torsion free subgroup of $\hbox{\rm Aut}(X)$. Then $\Gamma$ is mixing iff every element
$\gamma\in \Gamma-\{e\}$ is ergodic.
\end{pro}

\begin{proof}
If the action of $\Gamma$ on $X$ is mixing, then the action of every
infinite subgroup of $\Gamma$ is mixing as well, and in particular, every
$\gamma\in \Gamma-\{e\}$ is ergodic.

Conversely, suppose that the action of $\Gamma$ on $X$ is not mixing.
Then for some $(\chi,\psi)\in \hat X^2-\{(0,0)\}$,
the set
$$
S=\{\gamma\in\Gamma: \hat \gamma \chi=\psi\}
$$
is infinite. For every $\gamma\in S^{-1} S$, we have $\hat
\gamma\chi=\chi$, and the action of such $\gamma$ on $X$ is not ergodic.
This proves the proposition.
\end{proof}

Now we assume that {\it $X$ is connected and $\dim X=d<\infty$.}
We are going to show that under these assumptions, the torsion free condition in
Proposition \ref{p:mix} can be omitted. But first, we need the
following  lemma (see \cite[Lemma~4.3]{ber0} for a different proof).

\begin{lem}\label{l:finite}
Every torsion subgroup  (i.e., every element is of finite order) of $\hbox{\rm GL}(d,\mathbb{Q})$ is finite.
\end{lem}

In the proofs below, we use some basic facts about algebraic groups
and Zariski topology, which can be found in \cite{mor} and \cite{sp}.

\begin{proof}
Let $\Gamma$ be a torsion subgroup of $\hbox{\rm GL}(d,\mathbb{Q})$. The eigenvalues of a matrix in $\Gamma$
are roots of unity each having degree at most $d$ over $\mathbb{Q}$. Hence, their order is bounded,
and there exists $n\ge 1$ such that $\Gamma^n=\{e\}$. 
Let $G\subset \hbox{SL}(d,\mathbb{C})$ be the Zariski closure of $\Gamma$. Then its connected component $G^o$ has finite index in $G$,
and $G^n=\{e\}$. For $g\in G$, let $g=g_sg_u$ be the Jordan decomposition of $g$. Since $g_u$ is unipotent and $g_u^n=e$,
it follows that $g_u=e$ and every element of $G$ is semisimple. Hence, $G^o$ is a torus and since $(G^o)^n=\{e\}$,
we deduce that $G^o=\{e\}$ and $\Gamma$ is finite.
\end{proof}

\begin{pro}\label{c:mix}
Let $X$ be a compact connected finite-dimensional abelian group and
$\Gamma$ an infinite subgroup of $\hbox{\rm Aut}(X)$. Then the
following statements are equivalent:
\begin{enumerate}
\item[(a)] The action of $\Gamma$ on $X$ is mixing.
\item[(b)] Every infinite cyclic subgroup of $\Gamma$ is ergodic on $X$.
\item[(c)] For every $\gamma\in\Gamma$ of infinite order, the linear map $\hat
  \gamma$ does not have roots of unity as eigenvalues.
\end{enumerate}
\end{pro}

\begin{proof}
It is well-known that (a)$\Rightarrow$(b) and (b)$\Leftrightarrow$(c).
To show that (b)$\Rightarrow$(a), we observe that if the action of
$\Gamma$ on $X$ is not mixing, then for some $\chi\in \hat X-\{0\}$,
the subgroup $\{\gamma\in\Gamma:\, \hat\gamma\chi=\chi\}$ is infinite
(see the proof of Proposition \ref{p:mix}), and it suffices to show that 
this subgroup contains an element of infinite order.
This follows from Lemma \ref{l:finite}.
\end{proof}

Note that Proposition \ref{c:mix}((a)$\Leftrightarrow$(b)) fails in 
general if $X$ is disconnected or infinite-dimensional (see Example
\ref{ex:1} below). Also, it fails for semigroups (see Example
\ref{ex:2} below).

The following lemma is used in the proof of Corollary \ref{p:bhat}.

\begin{lem}\label{c:solv}
Every solvable mixing subgroup of $\hbox{\rm Aut}(X)$ is a finite extension of abelian group.  
\end{lem}

\begin{proof}
Let $\Gamma$ be a solvable mixing subgroup of $\hbox{\rm Aut}(X)$. We
show  that
$\hat\Gamma\subset\hbox{GL}(d,\mathbb{Q})$ is a finite extension of abelian group.
Let $G\subset \hbox{GL}(d,\mathbb{C})$ be the Zariski closure of $\Gamma$.
Then $G$ is solvable too. The connected component ${G}^o$ has finite index in $G$, and 
is is conjugate to a subgroup of the upper triangular subgroup (see \cite[Section~6.3]{sp}). 
In particular, the commutant $[G^o,G^o]$ is a unipotent subgroup.
The subgroup  $\hat \Gamma_0={G}^o\cap \hat\Gamma$ has finite index in $\hat\Gamma$. Since $\Gamma$ is mixing, it follows from
Proposition \ref{c:mix} that $[\hat\Gamma_0,\hat\Gamma_0]=1$.
This proves the corollary.
\end{proof}

\begin{proof}[Proof of Corollary \ref{p:bhat}]
Note that the subgroup $\Gamma$ is isomorphic to a subgroup of $\hbox{GL}(d,\mathbb{Q})$.
By the Tits alternative (see \cite{tits} or \cite[Section~5J]{mor}), $\Gamma$ is either finite extension of solvable group or contains a nonabelean free subgroup.
Thus, we may assume that $\Gamma$ contains a nonabelean free subroup. Let $\gamma$ and $\delta$ be free generators
and let $T_i=\delta^{-i}\gamma\delta^i$. Then 
$$
T_i^nT_j^{-n}=\delta^{-i}\gamma^n\delta^{i-j}\gamma^{-n}\delta^j\to
\infty\quad\hbox{for $i\ne j$.}
$$
On the other hand, linear maps $\hat T_i$ have the same characteristic polynomial. Hence, it follows from
Corollary \ref{th:cond_gen}(b) (or Proposition \ref{p_not_mix})
that the set $\{T_1,\ldots, T_{d+1}\}$ is not mixing.
This implies that $\Gamma$ is not mixing of order $d+1$.
\end{proof}

Using Proposition \ref{c:mix}, we develop two approaches to construction of mixing subgroups.
The first approach is based on the result of Y. Benoist \cite{ben} on
asymptotic cones of discrete groups 
(see Proposition \ref{p:mix_even})
and
the second approach is based on the theory of division algebras (see
Corollary \ref{c:mix_square}).  

\begin{pro}\label{p:mix_even}
For every even $d\ge 2$, there exists a mixing subgroup of $\hbox{\rm Aut}(\mathbb{T}^d)$
which is Zariski dense in $\hbox{\rm SL}(d,\mathbb{C})$. 
\end{pro}

\begin{proof}
We start by reviewing a result of Y.~Benoist from \cite{ben}, which
will be used in the proof.

For $g\in\hbox{SL}(d,\mathbb{R})$, let us denote by
$\lambda_1(g),\ldots,\lambda_d(g)$ the eigenvalues of $g$ such that $|\lambda_1(g)|\ge\cdots\ge |\lambda_d(g)|$
and
$$
\ell_g=(\log|\lambda_1(g)|,\ldots,\log |\lambda_d(g)|).
$$
The vector $\ell_g$ belongs to the set
$$
\mathfrak{a}^+\stackrel{def}{=}
\left\{(x_1,\ldots,x_d)\in\mathbb{R}^d:\,\sum_{i=1}^dx_i=0,\,x_1\ge \cdots \ge x_d\right\}.
$$
Let $\Gamma$ be a subgroup of $\hbox{SL}(d,\mathbb{R})$.
The {\it limit cone} $\ell_\Gamma$ of $\Gamma$ is the smallest closed cone in $\mathfrak{a}^+$ that
contains all $\ell_\gamma$, $\gamma\in\Gamma$. Since $\Gamma$ is a group, the limit cone $\ell_\Gamma$
is stable under the involution
$$
i(x_1,\ldots,x_d)=(-x_d,\ldots,-x_1).
$$
It was shown by Y. Benoist in \cite{ben} that if $\Gamma$ is Zariski dense, then the asymptotic cone $\ell_\Gamma$ is convex, has nonempty interior, and
is equal to the {\it asymptotic cone} of $\Gamma$. The asymptotic cone is the cone consisting of limit directions
of the set 
$$\{\log(\mu(\gamma)):\, \gamma\in\Gamma\}\subset \mathfrak{a}^+
$$
 where $\mu(g)$ denotes the $A^+$ component of $g$
with respect to $KA^+K$-decomposition ($K=\hbox{SO}(n)$,
 $A^+=$ positive Weyl chamber).
In the case when $\Gamma$ is a lattice, the asymptotic cone is always equal to $\mathfrak{a}^+$.
In particular,
$$
\ell_{\hbox{\tiny SL}(d,\mathbb{Z})}=\mathfrak{a}^+.
$$
If $\Gamma$ is a Zariski dense subgroup, Y. Benoist  also showed in \cite{ben} that
for every closed convex $i$-invariant cone $\mathcal{C}\subset\ell_\Gamma$ with nonempty interior,
there exists a Zariski dense subgroup $\Gamma_0\subset \Gamma$ such that $\ell_{\Gamma_0}=\mathcal{C}$.

Suppose that $d=2k$. For $(x_1,\ldots,x_{2k})\in\mathfrak{a}^+$,
$$
x_k\le -k^{-1}\sum_{i=k+1}^{2k} x_i,\quad x_{k+1}\ge k^{-1}\sum_{i=1}^{k} x_i,
$$
and for any $\delta\in (0,k^{-1})$
$$
\mathcal{C}_\delta=\left\{(x_1,\ldots,x_{2k})\in\mathfrak{a}^+:\,\begin{tabular}{l}$x_k\ge 0\ge x_{k+1}$\\
$x_k\ge -\delta\sum_{i=k+1}^{2k} x_i$\\ $x_{k+1}\le -\delta\sum_{i=1}^{k} x_i$\end{tabular}\right\}
$$
is closed convex $i$-invariant cone with nonempty interior. Hence, there exists a Zariski dense subgroup
$\Gamma$ of $\hbox{SL}(d,\mathbb{Z})$ such that $\ell_\Gamma=\mathcal{C}_\delta$. Since
$$
\mathcal{C}_\delta\cap \{x_i=0\}=0
$$
for every $i=1,\ldots,2k$, the group $\Gamma$ contains no element (except identity) with an eigenvalue of absolute value one.
By Proposition \ref{c:mix}, $\Gamma$ is mixing. This proves the proposition.
\end{proof}

If $\Gamma$ is a finitely generated subgroup of the automorphism group of a compact connected finite-dimensional
abelian group $X$, then by the Selberg lemma
(see \cite[Theorem~4.1]{cass} or \cite[Section 5I]{mor}), $\Gamma$ contains a torsion free subgroup of finite index.
Clearly, this subgroup is mixing iff $\Gamma$ is mixing. For torsion free subgroup Proposition
\ref{c:mix} can be restated as follows:

\begin{pro}\label{p:2}
Let $\Gamma$ be a torsion free subgroup of
$\hbox{\rm Aut}(X)$. Then the action of $\Gamma$ on $X$ is mixing iff
$$
\hat\Gamma-\hat\Gamma\subset \{0\}\cup\hbox{\rm GL}(d,\mathbb{Q}).
$$
\end{pro}

Recall that the Jacobson radical of a ring (with a unit) $R$ is the
intersection of all maximal ideals of $R$. We denote by $R^\times$ the
group of units of a ring $R$.

Let $A_\Gamma\subset \hbox{M}(d,\mathbb{Q})$ be the $\mathbb{Q}$-span of $\hat\Gamma$,
$J_\Gamma\subset A_\Gamma$ the Jacobson radical of $A_\Gamma$ and
$$
\pi:A_\Gamma\to A_\Gamma/J_\Gamma
$$
the factor map.

\begin{pro}\label{p:jacobson}
Let $\Gamma$ be a torsion free subgroup of $\Aut(X)$. Then the action of $\Gamma$ on $X$ is mixing iff
$$
\hat\Gamma\cap (1+J_\Gamma)=1\quad\hbox{and}\quad \pi(\hat\Gamma)-\pi(\hat\Gamma)\subset \{0\}\cup(A_\Gamma/J_\Gamma)^\times.
$$
\end{pro}

\begin{proof}
Recall that the Jacobson radial is nilpotent and $1+J_\Gamma\subset A_\Gamma^\times$.

Suppose that action of $\Gamma$ on $X$ is mixing. Since $1+J_\Gamma$ consists of unipotent matrices,
it follows from Proposition \ref{c:mix} that $\hat\Gamma\cap (1+J_\Gamma)=1$. The second property follows
Proposition \ref{p:2}.

Conversely, suppose that these properties are satisfied. If for some $a\in A_\Gamma$, $\pi(a)$ is invertible,
then there exists $b\in A_\Gamma$ such that $ab\in 1+J_\Gamma$ and it follows that $a$ is invertible as well.
Therefore,
$$
\hat\Gamma-\hat\Gamma\subset J_\Gamma\cup A_\Gamma^\times
$$
If $\gamma_1-\gamma_2\in J_\Gamma$ for some $\gamma_1,\gamma_2\in\hat\Gamma$, then $\gamma_1^{-1}\gamma_2\in 1+J_\Gamma$
and $\gamma_1=\gamma_2$. This shows that 
$$
\hat\Gamma-\hat\Gamma\subset \{0\}\cup \hbox{GL}(d,\mathbb{Q})
$$
and the action of $\Gamma$ on $X$ is mixing by Proposition \ref{p:2}. 
\end{proof}

Proposition \ref{p:jacobson} implies, in particular, that the action of $\Gamma$ on $X$ is mixing provided that the $\mathbb{Q}$-span
of $\hat\Gamma$ is a division subalgebra $D$ in $\hbox{M}(d,\mathbb{Q})$. This is possible only when $d=(\dim D)^l$
for some $l\ge 1$. In particular, $d$ must be a perfect square.

\begin{cor}\label{c:mix_square}
For every perfect square  $d>1$, there exists a mixing subgroup $\Gamma$ of $\hbox{\rm Aut}(\mathbb{T}^{d})$
such that the  Zariski closure of $\Gamma$ is conjugate to
$$
\{\hbox{\rm diag}(g,\ldots,g):\, g\in\hbox{\rm SL}(\sqrt{d},\mathbb{C})\}.
$$
The subgroup $\Gamma$ is not mixing of order $p+1$ where $p$ is the smallest prime divisor of $\sqrt{d}$. 
\end{cor}

\begin{proof}
There exists a central division algebra $D$ over $\mathbb{Q}$ 
such that $\dim_\mathbb{Q} D=d$ and $D$ is split over $\mathbb{R}$.
Denote by $\hbox{SL}(1,D)$ the group consisting of elements of $D$ whose
reduced norm is equal to one.
Consider the right and left regular representations
\begin{align*}
\rho:&D\otimes\mathbb{C}\to\hbox{End}(D\otimes\mathbb{C}):\; \rho(d)x=x\cdot d,\\
\lambda:&D\otimes\mathbb{C}\to\hbox{End}(D\otimes\mathbb{C}):\; \lambda(d)x=d\cdot x.
\end{align*}
Let $\mathcal{O}$ be an order in $D$.
Note that $(D\otimes\mathbb{R})/\mathcal{O}$ can be identified with the torus $\mathbb{T}^d$ and with respect to a basis of $\mathcal{O}$, we have 
\begin{align*}
G&\stackrel{def}{=}\rho(\hbox{SL}(1,D\otimes\mathbb{R}))\subset \hbox{SL}(d,\mathbb{R}),\\
\Gamma&\stackrel{def}{=}\rho(\hbox{SL}(1,\mathcal{O}))\subset \hbox{SL}(d,\mathbb{Z}).
\end{align*}

Since $D$ splits over $\mathbb{R}$,  $G\simeq \hbox{SL}(k,\mathbb{R})$ where $d=k^2$.
By the Borel--Harish--Chandra theorem (see \cite[Ch. IV]{pr}), $\Gamma$ is a lattice in $G$. 
This implies that the Zariski closure of $\Gamma$ is $\rho(\hbox{SL}(1,D\otimes\mathbb{C}))$. 
Note that 
$$
D\otimes \mathbb{C}\simeq \hbox{M}(k,\mathbb{C}),\quad \hbox{SL}(1,D\otimes\mathbb{C})\simeq \hbox{SL}(k,\mathbb{C}),
$$
and as a $\rho(D\otimes \mathbb{C})$-module,
$$
D\otimes\mathbb{C}\simeq \mathbb{C}^k\oplus\cdots\oplus \mathbb{C}^k.
$$
Hence,
$$
\rho(\hbox{SL}(1,D\otimes\mathbb{C}))=\{\hbox{diag}(g,\ldots,g):\, g\in\hbox{SL}(k,\mathbb{C})\}
$$
in a suitable basis.

Let $\hbox{Tr}:D\to \mathbb{Q}$ denote the reduced trace
of the division algebra $D$, and 
$\hat{\mathcal{O}}\subset D$ the dual order of $\mathcal{O}$, that is,
$$
\hat{\mathcal{O}}=\{u\in D_\mathbb{Q}:\, \hbox{Tr}(\mathcal{O}\cdot u)\in\mathbb{Z}\}.
$$
Then the set of characters of $D/\mathcal{O}$ is indexed by $\hat{\mathcal{O}}$:
$$
\chi_u(x)=\exp( 2\pi i\,\hbox{Tr}(x\cdot u)),\quad u\in\hat{\mathcal{O}},
$$
and the dual action of $\Gamma$ is
$$
u\mapsto \gamma\cdot u, \quad u\in\hat{\mathcal{O}},\; \gamma\in\Gamma.
$$
Hence, the $\mathbb{Q}$-span of $\hat\Gamma$ is equal to the $\mathbb{Q}$-span of
$\lambda(\hbox{SL}(1,\mathcal{O}))$, and since $\hbox{SL}(1,\mathcal{O})$ is Zariski dense in $\hbox{SL}(1,D\otimes\mathbb{C})$,
it is equal to $\lambda(D)$. Now it follows from Proposition \ref{p:jacobson} that the (right) action of $\Gamma$ on $D/\mathcal{O}$ is mixing.

The central division algebra $D$ contains a splitting field $F$ such that $F/\mathbb{Q}$
is a cyclic extension of degree $k$. Moreover, since $D$ splits over $\mathbb{R}$, $F$ can be taken to be real.
Then $F$ contains a Galois subfield $E$ such that $|E:\mathbb{Q}|=p$.
By Dirichlet theorem (see \cite[Ch. 2]{koch}), $E$ contains a unit $\gamma$ of infinite order (unless $E/\mathbb{Q}$ is a complex quadratic
extension, which is not the case). Since
$$
\hbox{N}(\gamma)=N_{E/\mathbb{Q}}(\gamma)^{\dim_E D} =\pm 1,
$$
we may choose $\gamma\in \Gamma$. There exist $\alpha_i\in\mathbb{Q}$, $i=0,\ldots,p$, $\alpha_p=1$,
such that
$$
\sum_{i=0}^p \alpha_i \gamma^{i}=0.
$$
We claim that there exists $\gamma_n\in \Gamma$ such that
$$
\gamma_n^{-1}\gamma^i\gamma_n\to\infty\quad\hbox{ as $n\to\infty$.}
$$
for $i=1,\ldots,p$. It suffices to check that the centralizer $C_\Gamma (\gamma^{p!})$ has infinite index in $\Gamma$.
If this is not the case, then  $C_\Gamma (\gamma^{p!})$ is a lattice in $G$, and it follows that $\gamma^{p!}$ 
lies in the center of $D$, which is a contradiction. We have
$$
\sum_{i=0}^p  (\gamma_n^{-1}\gamma^{i}\gamma_n)\alpha_i=0.
$$
Since $\ell\alpha_i\in\hat{\mathcal{O}}$ for some $\ell\in \mathbb{N}$, this proves that $\rho(\Gamma)$ is not mixing of order $p+1$.
\end{proof}

\begin{cor}\label{c:mix_new}
For every $d\ge 2$, $d\ne 3,5,7,9$, there exists a free nonabelian mixing subgroup of $\hbox{\rm Aut}(\mathbb{T}^d)$
which is not mixing of order 3.
\end{cor}

\begin{proof}
Let $d=d_1+d_2$ where $d_1>1$ is a perfect square and $d_2>1$ is even. If $d\ne 9$, we may take $d_2\ge 2$.
By Proposition \ref{p:mix_even} and Corollary \ref{c:mix_square}, $\hbox{Aut}(\mathbb{T}^{d_i})$, $i=1,2$,
contain not virtually abelian free mixing subgroups.
Hence, by the Tits alternative, there exist injective homomorphisms 
$$\phi_i:F_2\to\hbox{Aut}(\mathbb{T}^{d_i}),\quad i=1,2,$$ where $F_2$
denotes the free group with 2 generators such that the action of $\phi_i(F_2)$ on $\mathbb{T}^{d_i}$ is mixing. 
Let
$$
\Gamma=\{(\phi_1(\delta),\phi_2(\delta)):\gamma\in F_2\}\subset\hbox{Aut}(\mathbb{T}^{d_1}\times\mathbb{T}^{d_2})=\hbox{Aut}(\mathbb{T}^{d}).
$$
If $\Gamma$ is not mixing, there exist $x,y\in\mathbb{Q}^d-\{0\}$ and $\gamma_n\in\Gamma$ such that $\gamma_n\to\infty$
and ${}^t\gamma_n x=y$. Write $x=x_1+x_2$ and $y=y_1+y_2$ with respect to the decomposition $\mathbb{Q}^d=\mathbb{Q}^{d_1}\oplus\mathbb{Q}^{d_2}$
Then for some $i=1,2$, we have $x_i,y_i\in\mathbb{Q}^{d_i}-\{0\}$ and ${}^t\phi_i(\delta_n)x_i=y_i$ where $\delta_n\in F_2$ corresponds
to $\gamma_n$. This is a contradiction since the subgroups
$\phi_i(F_2)$ are mixing. It follows that $\Gamma$ is mixing.

If $d\ne 9$, the group $\Gamma$ preserves the direct product decomposition $\mathbb{T}^d=\mathbb{T}^{d-2}\times \mathbb{T}^2$.
If the action $\Gamma$ on $\mathbb{T}^d$ is mixing of order 3, then the restriction of this action to $\mathbb{T}^2$ is
mixing as well. However, this contradicts Proposition \ref{p:bhat}.
\end{proof}

\begin{cor}
There exists a not virtually abelian mixing subgroup of $\hbox{\rm Aut}(\mathbb{T}^9)$
which is not mixing of order 3.
\end{cor}

\begin{proof}
Let $K=\mathbb{Q}(\alpha)$ with $\alpha=\zeta+\bar\zeta$ where $\zeta$ is a primitive root of unity
of order $7$,  and $D\supset K$ be a central division algebra over $\mathbb{Q}$ with $\dim_\mathbb{Q} D=9$
(such algebra can be constructed using the cross product construction).
One can check that  $\alpha$ is a root of $x^3+x^2-2x-1=0$. In particular, this implies that 
$\alpha$ and $\beta=-1-\alpha$ are units in $K$ and $\alpha,\beta\in\hbox{SL}(1,D)$.
Let $\mathcal{O}$ be an order in $D$ that contains $\alpha$ and $\beta$. 
Using the argument from the proof of Corollary \ref{c:mix_square}, one can find a sequence
$\{\gamma_n\}\subset\hbox{SL}(1,\mathcal{O})$ such that $\gamma_n^{-1}\alpha\gamma_n\to\infty$ as $n\to\infty$.
Then
$$
\gamma_n^{-1}\alpha\gamma_n+\gamma_n^{-1}\beta\gamma_n+1=0
$$
and $\gamma_n^{-1}\beta\gamma_n\to\infty$, $\gamma_n^{-1}(\alpha\beta^{-1})\gamma_n\to\infty$.
As in the proof of Corollary \ref{c:mix_square}, this implies that the action of $\hbox{SL}(1,\mathcal{O})$
on $D/\mathcal{O}$ by right multiplication is not mixing of order 3.
\end{proof}

\begin{que}
Does there exist a not virtually abelian mixing subgroup in $\hbox{\rm
  Aut}(\mathbb{T}^d)$ for $d=3,5,7$?
\end{que}

According to Corollary \ref{p:bhat}, $\hbox{Aut}(\mathbb{T}^2)$
contains no free nonabelian subgroup which is mixing of order 3
(see also Proposition~2.31 in the electronic version of \cite{bg}).

\begin{que}
Is there a free nonabelian mixing of order 3 subgroup in $\hbox{\rm Aut}(\mathbb{T}^d)$ for  some $d\ge 3$?
\end{que}

Note that there exist free nonabelian semigroups of epimorphisms of the torus
$\mathbb{T}^d$ which are mixing of all orders (see Example
\ref{ex:mix2} below).



\section{Ergodicity and mixing}\label{sec:erg}

In this section, we prove Theorems \ref{th:solv1}, \ref{th:solv2}, and
\ref{th:tits}. 

First, we recall the following well-known
characterization of ergodicity (see, for example, \cite[Chapter I]{sch}):

\begin{pro}\label{p:ergodic}
Let $\Gamma$ be a groups of automorphisms of a compact abelian group
$X$. Then the action of $\Gamma$ on $X$ is ergodic iff the action of
$\Gamma$ on $\hat X$ has no finite orbits except the trivial character.
\end{pro}

Using Proposition \ref{p:ergodic}, we deduce

\begin{pro}\label{l:irr}
Let $X$ be a compact connected abelian group and
$\Gamma\subset\hbox{\rm Aut}(X)$. Then if the action of $\Gamma$ on $X$
is strongly irreducible, then it is ergodic.
\end{pro}

Note that the converse of Proposition \ref{l:irr} is not true (see Example
\ref{ex:not_irr} below).

\begin{proof}
Suppose that the action of $\Gamma$ on $X$ is not ergodic. Then there
exist $\chi\in\hat X-\{0\}$ and a subgroup $\Lambda$ of finite index
in $\Gamma$ such that $\Lambda\cdot\chi=\chi$. Consider the subgroup
$$
A=\{\psi\in \hat X:\, k\psi\in\mathbb{Z}\cdot\chi\quad\hbox{for some $k\ge
  1$}\}.
$$
Using that $\hat X$ is torsion free, one can check that $\Lambda$ acts
trivially on $A$. In particular, $A\ne \hat X$.
Also, it is clear that $\hat X/A$ is torsion free.
Hence, there exists a proper closed connected
$\Lambda$-invariant subgroup
$$
\{x\in X:\, \chi(x)=1\quad\hbox{for all $\chi\in A$}\},
$$
and the action of $\Gamma$ on $X$ is not
strongly irreducible.
\end{proof}

In the proofs of the Theorems \ref{th:solv1}, \ref{th:solv2}, and
\ref{th:tits}, we will need the following three lemmas.

\begin{lem}\label{l:hom}
Let $\Gamma$ be a group and $\rho_i:\Gamma\to (\mathbb{C},+)$,
 $i=1,\ldots,t$, nontrivial homomorphisms. Then there exists
 $\gamma\in \Gamma$ such that $\rho_i(\gamma)\ne 0$ for every
 $i=1,\ldots,t$.
Moreover, the set
$$
R=\{\gamma\in\Gamma:\,\hbox{$\rho_i(\gamma)\ne 0$ for every
 $i=1,\ldots,t$}\}.
$$
generates $\Gamma$.
\end{lem}

\begin{proof}
Consider a homomorphisms $\rho:\Gamma\to \mathbb{C}^t$ defined by
$$
\rho(\gamma)=(\rho_1(\gamma),\ldots,\rho_t(\gamma)).
$$
Then $\Delta=\rho(\Gamma)$ is a subgroup of $\mathbb{C}^t$ such
that $\pi_i(\Delta)\ne 0$, $i=1,\ldots,t$,
where $\pi_i:\mathbb{C}^t\to \mathbb{C}$ the coordinate projection.
It suffices to show that
$$
\Delta\nsubseteq \bigcup_{i=1}^t \pi_i^{-1}(0).
$$
Suppose that this is not the case. Then the Zariski closure $\bar\Delta$
of $\Delta$ is a linear subspace of $\mathbb{C}^t$ and
$$
\bar \Delta= \bigcup_{i=1}^t (\bar \Delta\cap \pi_i^{-1}(0)).
$$
However, this equality is impossible because $\bar \Delta\cap \pi_i^{-1}(0)$
are proper linear subspaces of $\bar \Delta$. This contradiction
proves the first part of the lemma.

To prove the second part, take any $\gamma\in\Gamma$ and
$\delta\in R$. Then for $k\ge 1$,
$$
\rho(\gamma\delta^k)=\rho(\gamma)+k\rho(\delta),
$$
and taking $k$ such that 
$$
k\ne -\pi_i(\rho(\gamma))/\pi_i(\rho(\delta))\quad \hbox{ for every
 $i=1,\ldots,t$},
$$
we have $\delta^k\in\Gamma_0$ and $\gamma\delta^k\in R$.
Hence, $\gamma\in\left<R\right>$. This proves the lemma.
\end{proof}

\begin{lem}\label{l:chain}
Let $\Gamma$ be a solvable subgroup of $\hbox{\rm GL}(d,\mathbb{Q})$.
Then there exist a subgroup $\Lambda$ such that $|\Gamma:\Lambda|<\infty$
and the commutant $\Lambda'$ is unipotent, and a flag 
$$
\mathbb{Q}^d=V_1\supset V_2\supset \cdots \supset V_{s+1}=\{0\}
$$
consisting of rational $\Lambda$-invariant subspaces  such that
$\Lambda|_{V_i/V_{i+1}}$ is abelian for all $i=1,\ldots,s$. 
\end{lem}

\begin{proof}
There exists a subgroup $\Lambda$ of finite index in $\Gamma$ which
can be conjugated (over $\mathbb{C}$) to a subgroup of the group of the upper
triangular matrices  (see, for example, the proof of Lemma \ref{c:solv}). Then the commutant $\Lambda'$ is a
unipotent subgroup. Hence, the subspace $V^{\Lambda'}$
of $\Lambda'$-invariant vectors is not trivial.  Since
$\Lambda'$ is normal in $\Lambda$, this subspace is
$\Lambda$-invariant. Also, it is clear that $V^{\Lambda'}$  is rational,
and $\Gamma|_{V^{\Lambda'}}$ is abelian. Now the lemma
follows by induction on dimension.
\end{proof}

For a subgroup $\Gamma\subset\hbox{GL}(d,\mathbb{Q})$, we denote by
$\bar\Gamma$ its Zariski closure and by $\bar\Gamma^\circ$ the
connected component of the closure.

\begin{lem}\label{l:fg}
Every subgroup $\Gamma$ of $\hbox{\rm GL}(d,\mathbb{Q})$ contains a
finitely generated subgroup $\Lambda$ such that $\bar\Lambda=\bar\Gamma$.
\end{lem}

\begin{proof}
Take a finitely generated subgroup $\Delta$ such that
$\dim \bar\Delta^\circ$ is maximal among all finitely generated
subgroups. Then for every $\gamma\in \Gamma$, 
$$
\overline{\left<\Delta,\gamma\right>}^o=\bar\Delta^o.
$$
In particular, $\gamma^{-1}\bar\Delta^o\gamma\subset \bar\Delta^o$,
and the group $\Gamma\cap \bar\Delta^o$ is normal in $\Gamma$.
Also, since $\overline{\left<\Delta,\gamma\right>}$ has finitely many
connected components, $\gamma^k\in \bar\Delta^o$ for some $k\ge 1$ and the group
$\Gamma/(\Gamma\cap \bar\Delta^o)$ consists of elements of finite order.
The algebraic group $\bar\Gamma/\bar\Delta^o$ is defined over $\mathbb{Q}$
and it embeds via a $\mathbb{Q}$-map into $\hbox{GL}(n)$ for
some $n\ge 2$. Under this map, the subgroup $\Gamma/(\Gamma\cap
\bar\Delta^o)$ is embedded into $\hbox{GL}(n,\mathbb{Q})$.
Hence, it follows from Lemma \ref{l:finite} that $\Gamma/(\Gamma\cap
\bar\Delta^o)$ is finite. This implies that $\bar\Delta$ has finite
index in $\bar\Gamma$. Since $\Gamma$ is dense in $\bar \Gamma$, every
coset of $\bar \Delta$ in $\bar\Gamma$ contains a representative from
$\Gamma$. Now the required group $\Lambda$ can be taken to be generated by
$\Delta$ and these coset representatives.
\end{proof}

\begin{proof}[Proof of Theorem \ref{th:solv1}]
(b)$\Rightarrow$(a): 
Suppose that there exists a closed connected virtually
$\Gamma$-invariant subgroup $Y$ of $X$ such that the action of
$\Gamma_Y$ on $Y$ is not ergodic. Then by Proposition \ref{p:ergodic},
there exists a subgroup $\Lambda$ with $|\Gamma:\Lambda|<\infty$ and 
$\chi\in \hat Y-\{0\}$ such that $\Lambda \chi=\chi$. The character
group $\hat Y$ is equal to $\hat X/A(Y)$ where 
$$
A(Y)=\{\chi\in \hat X:\, \chi(Y)=1 \}.
$$
Since $Y$ is connected, $\hat Y$ is torsion free, and it follows that
$\hat X/A(Y)$ embeds in $(\hat X\otimes \mathbb{Q})/(A(Y)\otimes
\mathbb{Q})$. Therefore, the character $\chi$ gives a nonzero vector
in $(\hat X\otimes \mathbb{Q})/(A(Y)\otimes \mathbb{Q})$ which is fixed by $\Lambda$.
This implies that every element of $\Lambda$ has eigenvalue one, and 
by Proposition \ref{c:mix}, $\Lambda$ contains no mixing subgroup.
Since $\Lambda$ has finite index in $\Gamma$, $\Gamma$ does not
contain any mixing subgroups as well.

(b)$\Rightarrow$(a): 
First, we can pass to a finite index subgroup $\Lambda$ of $\Gamma$ as in Lemma
\ref{l:chain}. Then for every $i=1,\ldots, s$,
$$
V_i/V_{i+1}=\bigoplus_{j=1}^{n_i} (V_i/V_{i+1})_{\alpha_{i,j}},
$$
where $(V_i/V_{i+1})_{\alpha_{i,j}}$ denotes the weight space corresponding
to a homomorphism $\alpha_{i,j}:\Lambda\to \mathbb{C}^\times$.

Suppose that for some $\alpha_{i,j}$, the set
$\alpha_{i,j}(\Lambda)$ consists of roots of unity. Since $\alpha_{i,j}(\Lambda)$
consists of eigenvalues of matrices in $\hbox{GL}(d,\mathbb{Q})$,
it follows that for every $\alpha\in\alpha_{i,j}(\Lambda)$,
$[\mathbb{Q}(\alpha):\mathbb{Q}]\le d$ and $\alpha^N=1$ where $N\ge
1$ depends only on $d$.
Hence, passing again, if necessary, to a finite index subgroup if needed, we can assume that
$\alpha_{i,j}(\Lambda)=\{1\}$. Then there exists $v\in V_i(\mathbb{Q})-V_{i+1}(\mathbb{Q})$ such that
$$
\Lambda\cdot v=v+V_{i+1}.
$$
Let 
$$
Y=\{x\in X:\, \chi(x)=0\quad\hbox{for all $\chi\in  \hat X\cap V_{i+1}$}\}.
$$
$Y$ is a closed subgroup of $X$ with the character group equal to $\hat
X/(\hat X\cap V_{i+1})$. Since the character group of $Y$ is torsion free,
the group $Y$ is connected. Take $l\ge 1$ such that $lv\in \hat X$. This gives a nontrivial
character of $Y$ which is fixed by $\Lambda$. Hence, the action of $\Lambda$ on
$Y$ is not ergodic, which contradicts hereditary ergodicity.

It follows that there exists a finitely generated subgroup $\Lambda_0$
of $\Lambda$ such that for every $\alpha_{i,j}$, the set
$\alpha_{i,j}(\Lambda_0)$ contains an element which is not a root of
unity. Denote by $K$ the field generated  by the sets
$\alpha_{i,j}(\Lambda_0)$, $i=1,\ldots,s$, $j=1,\ldots,n_i$. Since $\Lambda_0$ is finitely generated,
$[K:\mathbb{Q}]<\infty$. By Kronecker's lemma (Lemma \ref{l:abs_1}),
for every $\alpha_{i,j}$ there exists a an absolute value  $|\cdot|_{i,j}$ of
the field $K$ such that $|\alpha_{i,j}(\Lambda_0)|_{i,j}\ne 1$.
Consider the set of nontrivial homomorphisms
\begin{equation}\label{eq:rho}
\rho_{i,j}(\lambda)=|\alpha_{i,j}(\lambda)|_{i,j}: \Lambda_0\to
\mathbb{R}^+,\quad i=1,\ldots,s,\, j=1,\ldots, n_i.
\end{equation}
By Lemma \ref{l:hom}, there exists $\gamma\in \Lambda_0$ such that
$\rho_{i,j}(\gamma)\ne 1$ for all $\rho_{i,j}$'s.
In particular, $\gamma$ has no roots of unity as eigenvalues.
Hence, it is ergodic, and moreover, it is mixing of all orders by
Rokhlin's theorem (Theorem \ref{th:roh}).
\end{proof}

\begin{proof}[Proof of Theorem \ref{th:solv2}]
Note that (b)$\Rightarrow$(a) follows from Theorem \ref{th:solv1}, and
it suffices to prove that (a)$\Rightarrow$(b).

By Lemma  \ref{l:fg}, there exists a finitely generated subgroup
$\Delta$ of $\Gamma$ such that $\bar \Delta=\bar\Gamma$. Since
$\Gamma$ is not virtually nilpotent, $\Delta$ is not virtually
nilpotent as well. By Theorem \ref{th:solv1}, 
$\Gamma$ contains a mixing transformation $\gamma_0$.
Then the group generated by $\Delta$ and $\gamma_0$ is finitely
generated, not virtually nilpotent, and it satisfies (a).
Hence, we can assume that $\Gamma$ is finitely generated.

Let $\Lambda$ be a finite index subgroup of $\Gamma$ as in Lemma
\ref{l:chain} and 
$$\alpha_{i,j}:\Gamma\to\mathbb{C}^\times,\;\;
i=1,\ldots,s,\; j=1,\ldots,n_i,
$$
the weights of the action of $\Lambda$
on $V_i/V_{i+1}$. Denote by $K$ the field generated by the sets
$\alpha_{i,j}(\Lambda)$,  $i=1,\ldots,s$, $j=1,\ldots,n_i$.
Since $\Lambda$ is finitely generated, $K$ has finite degree over $\mathbb{Q}$.
As in the proof of Theorem \ref{th:solv1}, we deduce from (a) that
for every $\alpha_{i,j}$ there exists a absolute value $|\cdot|_{i,j}$
of the field $K$ such that the homomorphism
$$
\rho_{i,j}(\lambda)=|\alpha_{i,j}(\lambda)|_{i,j}:\Lambda\to \mathbb{R}^+
$$
is not trivial. Set
$$
R=\{\lambda\in\Lambda:\, \rho_{i,j}(\lambda)\ne 0 \hbox{ for all
  $i=1,\ldots,s$ and $j=1,\ldots,n_i$}\}.
$$
By Lemma \ref{l:hom}, $R$ generates $\Lambda$.
Note that every $\lambda\in R$, $\hat\lambda$ does not have roots of unity as
eigenvalues, and by Rokhlin's Theorem (Theorem \ref{th:roh}), $\lambda$ is
mixing of all orders.

\medskip 
\noindent {\it Claim.} There there exist $\delta\in R$ and $\mu\in\Lambda'$
such that the semigroup $S=\left<\delta,\delta\mu\right>$ is free.
\medskip

Consider the derived series of $\Lambda$:
$$
\Lambda\supset \Lambda'\supset \Lambda^{(2)}\supset\cdots \supset\Lambda^{(k+1)}=\{e\}.
$$
Suppose that $\Lambda^{(i)}/\Lambda^{(i+1)}$ is finitely generated for
$i=0,\ldots, l-1$, but $\Lambda^{(l)}/\Lambda^{(l+1)}$ is not finitely
generated. Then $\Lambda/\Lambda^{(l)}$ is polycyclic, and in
particular, finitely presented. Applying \cite[Lemma~4.9]{ros},
we deduce that there exists a finite subset $T$ of $\Lambda^{(l)}/\Lambda^{(l+1)}$
such that $\Lambda^{(l)}/\Lambda^{(l+1)}$ is generated by
$\lambda T\lambda^{-1}$, $\lambda\in\Lambda/\Lambda^{(l)}$.
Since $\Lambda/\Lambda^{(l)}$ is polycyclic, $\Lambda'/\Lambda^{(l)}$
is finitely generated. Also, $\Lambda'$ is nilpotent (see Lemma
\ref{l:chain}). This implies that the set 
$\lambda T\lambda^{-1}$, $\lambda\in\Lambda'/\Lambda^{(l)}$
generates a finitely generated subgroup of $\Lambda^{(l)}/\Lambda^{(l+1)}$.
Since $R$ generates $\Lambda$, there exist $\lambda_1,\ldots,\lambda_r\in R$ such that
$$
\Lambda=\lambda_1^\mathbb{Z}\cdots\lambda_r^\mathbb{Z}\Lambda'.
$$
It follows that there exists a finite set $Q\subset \Lambda^{(l)}/\Lambda^{(l+1)}$
such that the group $\Lambda^{(l)}/\Lambda^{(l+1)}$ is generated by
$$
\lambda_1^{n_1}\cdots\lambda_r^{n_r}q
\lambda_r^{-n_r}\cdots\lambda_1^{-n_1},\quad q\in Q,\;
n_1,\ldots,n_r\in\mathbb{Z}.
$$
Hence, since $\Lambda^{(l)}/\Lambda^{(l+1)}$ is not finitely
generated, we deduce that there exists $\delta\in R$ and $\mu\in
\Lambda^{(l)}/\Lambda^{(l+1)}$ such that
$$
\delta^{n}\mu\delta^{-n},\quad n\in\mathbb{Z},
$$
generates an infinitely generated subgroup. 
Now the claim follows from \cite[Lemma
4.8]{ros}.

Next, we consider the case when the all groups $\Lambda^{(i)}/\Lambda^{(i+1)}$ are finitely-generated.
Then
$$
\Lambda^{(i)}/\Lambda^{(i+1)}\simeq \mathbb{Z}^{d_i}\oplus A_i
$$
where $A_i$ is a finite abelian group. Denote by $\Delta_i$ the
preimage of $A_i$ under the factor map $\Lambda\to
\Lambda/\Lambda^{(i+1)}$. Note that $\Delta_i$ is a normal subgroup of
$\Lambda$. There exists a finite index subgroup $\Lambda_0$ of $\Lambda$
such that the action of $\Lambda_0$ on $\Lambda^{(i)}/\Delta_i\simeq
\mathbb{Z}^{d_i}$ is conjugate (over $\mathbb{C}$) to an action by
upper triangular matrices. For $\gamma\in\Lambda_0$, we denote by
$\beta_{i,j}(\gamma)$, $j=1,\ldots,n_i$, the eigenvalues of the
corresponding upper triangular matrix. Note that
the maps $\beta_{i,j}:\Lambda_0\to \mathbb{C}^\times$ are homomorphism.

Suppose that for every $i=1,\ldots,k$ and $j=1,\ldots,n_i$, the set
$\beta_{i,j}(\Lambda_0)$ consists of roots of unity. Since the sets
$\beta_{i,j}(\Lambda_0)$ consist of algebraic numbers of degree at
most $n_i$. It follows that there exists $N\ge 1$ such that for every
$\beta\in \beta_{i,j}(\Lambda_0)$, we have $\beta^N=1$. Hence, by passing
to a finite index subgroup, we may assume that
$\beta_{i,j}(\Lambda_0)=1$ for all $\beta_{i,j}$'s. Also, passing to a
finite index subgroup, we may assume
that $\Lambda_0$ acts trivially on $\Delta_i/\Lambda^{(i+1)}$.
Each of the linear maps
$$
\Lambda^{(i)}/\Delta_i\to\Lambda^{(i)}/\Delta_i: x\mapsto \gamma x\gamma^{-1},
\quad \gamma\in\Lambda_0, 
$$
is unipotent. This implies that the corresponding action $\Lambda_0$
is unipotent. Then this action is conjugate to the action by a
group of unipotent upper triangular matrices. Then the linear maps
$$
\Lambda^{(i)}/\Delta_i\to\Lambda^{(i)}/\Delta_i: x\mapsto \gamma x\gamma^{-1}x^{-1}=[\gamma,x],
\quad \gamma\in\Lambda_0, 
$$
generate a nilpotent subalgebra, and it follows that
$$
[\Lambda_0,\ldots,\Lambda_0,\Lambda^{(i)}]\subset \Delta_i\quad
\hbox{($n_i$ terms)}.
$$
Since $\Lambda_0$ acts trivially on $\Delta_i/\Lambda^{(i+1)}$, we
also have
$$
[\Lambda_0,\ldots,\Lambda_0,\Lambda^{(i)}]\subset \Lambda^{(i+1)}\quad
\hbox{($n_i+1$ terms)}.
$$
This implies that 
$$
[\Lambda_0,\ldots,\Lambda_0,\Lambda]=1,
$$
and in particular, $\Lambda_0$ is nilpotent, which is a contradiction.

We have shown that for some $i=1,\ldots,k$, $j=1,\ldots, n_i$, and
$\gamma\in\Lambda_0$, the number $\beta_{i,j}(\lambda)$ is not a root of unity.
Note that the numbers $\beta_{i,j}(\lambda)$, $j=1,\ldots, n_i$, are 
algebraic integers, and they are permuted by the action of the Galois
group. Hence, by Kronecker's lemma (Lemma \ref{l:abs_1}),
$|\beta_{i_0,j_0}(\lambda)|\ne 1$ for some $i_0=1,\ldots,k$ and
$j_0=1,\ldots, n_{i_0}$.
Note that since the action of $\lambda$ on
$\Lambda/\Lambda'$ is trivial, $i_0>1$.
By Lemma \ref{l:hom}, there exists $\lambda\in\Lambda_0$ such that
$|\beta_{i_0,j_0}(\lambda)|\ne 1$ and $\rho_{i,j}(\lambda)\ne 1$ for
all $\rho_{i,j}$'s as in (\ref{eq:rho}).

By \cite[Theorem~4.17]{ros}, there exists $\mu\in
\Lambda^{(i_0)}\subset \Lambda'$ such
that the semigroup
$S=\left<\lambda^n,\lambda^n\mu\right>$ is free for sufficiently large
$n\ge 1$. 
This proves the claim.

It remains to show that the action of the semigroup $S$ on $X$ is
mixing of all orders. Suppose that, in contrary, there
exist $x_j\in \hat X\otimes \mathbb{Q}$ and $\gamma_j^{(n)}\in \Gamma$,
$j=1,\ldots,t$, such that $\gamma_k^{(n)}\ldots\gamma_i^{(n)}\to
\infty$ for $1<i\le k\le t$ and
\begin{equation}\label{eq:mix_last}
\hat \gamma_1^{(n)} x_1+(\hat \gamma_1^{(n)}\hat \gamma_2^{(n)}) x_2+\cdots +(\hat \gamma_1^{(n)}\ldots\hat \gamma_t^{(n)})x_t=0.
\end{equation}
Denote by $p_i:V_j\to V_i/V_{i+1}$, $i=1,\ldots,s$, the projection maps. Since the
action of $\Gamma$ on $V_i/V_{i+1}$ is abelian  and $\mu\in \Lambda'$,
it follows that $\mu$ acts trivially on $V_i/V_{i+1}$. Hence, for
every $v\in V_i/V_{i+1}$,
$$
\hat \gamma_1^{(n)}\ldots\hat \gamma_k^{(n)}v=\delta^{l_1(n)+\cdots+l_k(n)}v.
$$
with $l_j(n)\to\infty$. Now we deduce from (\ref{eq:mix_last}) that
$$
\delta^{l_1(n)} p_1(x_1)+\delta^{l_1(n)+l_2(n)}
p_1(x_2)+\cdots +\delta^{l_1(n)+\cdots+l_t(n)} p_1(x_t)=0.
$$
According to our choice of $\delta$, the map $\delta$ has no roots of
unity as eigenvalues for the action on $V_1/V_2$. Therefore, it
follows from Rokhlin's Theorem (Theorem \ref{th:roh}) that
$p_1(x_j)=0$ and $x_j\in V_2$ for $j=1,\ldots,t$. Applying the same argument
to the spaces $V_i/V_{i+1}$ for $i=2,\ldots,s$, we deduce that $x_j=0$
for  every $j=1,\ldots,t$. This proves that the action of $S$ on $X$
is mixing of all orders.
\end{proof}

\begin{proof}[Proof of Theorem \ref{th:tits}]
Passing to a finite index subgroup, we may assume that the Zariski
closure $\bar \Gamma$ is connected.
By Lemma \ref{l:fg}, there exists a finitely generated subgroup
$\Lambda$ in $\Gamma$ such that $\bar\Lambda=\bar\Gamma$. 
In particular, $\Lambda$ is not virtually solvable.

Suppose that $\Lambda$ contains
a Zariski open solvable group $\Delta=\Lambda\cap U$, where $U$ is an
open subset of $\bar\Lambda$. Then for $\gamma_1,\gamma_2\in\Lambda$ such that
$\gamma_1^{-1}\gamma_2\notin \Delta$,
$$
\Lambda\cap \gamma_1U\cap \gamma_2U=\emptyset,
$$
and since $\Lambda$ is dense,
$$
\gamma_1U\cap \gamma_2U=\emptyset.
$$
Hence, we have a disjoint union
$$
\bar\Lambda=\bigcup_{\gamma\in \Gamma/\Delta} \gamma \Delta U.
$$
This implies that $\Delta=\Lambda$ and gives a contradiction.
Therefore, by \cite[Theorem~1.1]{br_g}, the group $\Lambda$ contains
nonabelian free subgroup $\Delta$ such that $\bar \Delta=\bar\Lambda=\bar\Gamma$.

Suppose that the action of $\Gamma$ on $X$ is ergodic, but the action
of $\Delta$ on $X$ is not ergodic. By Proposition
\ref{p:ergodic}, there exists $\chi\in \hat X-\{0\}$ such that
$\Delta\chi$ is finite. Then $\bar\Delta \chi=\bar\Gamma\chi$ is
finite, and this gives a contradiction.

For every closed connected subgroup $Y$ of $X$ and 
$$
A(Y)=\{\chi\in \hat X:\, \chi (Y)=1\},
$$
we have 
$$
A(Y)=\hat X\cap (A(Y)\otimes \mathbb{Q}).
$$
Since $\bar\Delta=\bar\Gamma$,  this implies that if $Y$ is
$\Delta$-invariant, then it is $\Gamma$-invariant.
In particular, this shows that if $\Gamma$ is strongly irreducible,
then $\Delta$ is strongly irreducible as well.

Suppose that the action of $\Delta$ on $X$ is not hereditarily
ergodic, i.e., there exist a closed connected virtually
$\Delta_0$-invariant subgroup $Y$, where $\Delta_0$ is a subgroup
of finite index in $\Delta$, and $\chi\in \hat Y-\{0\}$ such that
$\Delta_0\chi$ is finite. Then we deduce as above that $Y$
is invariant under $\Gamma_0=\Gamma\cap \bar \Delta_0$ which has
finite index in $\Gamma$. The character group of $Y$ can be identified
with $\hat X/A(Y)$. Moreover, since $Y$ is connected, $\hat X/A(Y)$ is
torsion-free, and the map
$$
i:\hat X/A(Y)\to (\hat X\otimes \mathbb{C}) /(A(Y)\otimes \mathbb{C})
$$
is injective. Using that $\Delta_0\chi$ is finite, we deduce that
$\Gamma_0\cdot i(\chi)$ and $\bar \Gamma_0\cdot i(\chi)$ are finite.
It follows that $\Gamma_0\chi$ is finite, and the
action of $\Gamma$ on $X$ is not hereditarily ergodic. This proves the theorem.
\end{proof}

\section{Examples}\label{sec:ex}

\begin{exa}[cf. Theorem \ref{th:general} and Corollary
  \ref{thm:limit}] \label{ex:last2}
For
$$
S=\left(\begin{matrix} 0 & -1\\ 1 & 0\end{matrix}\right)\quad
\hbox{and}\quad
T=\left(\begin{matrix} 0 & -1\\ 1 & -1\end{matrix}\right),
$$
we have 
\begin{enumerate}
\item[(a)] $S^4=id$ and $T^3=id$. In particular, the set $\{T,S\}$ is
  not mixing on $\mathbb{T}^2$.
\item[(b)] For every $l\ge 1$ such that
  $\phi(l)<(\dim\mathbb{T}^2)^2$,
the linear maps $\hat S^l$ and $\hat T^l$ don't have a
  common eigenvalue.
\item[(c)] There exists $f\in L^\infty(\mathbb{T}^2)$ such that the
  limit
$$
\lim_{n\to\infty} \int_{\mathbb{T}^2}f(S^{ln}x)f(T^{ln}x)\, dm(x)
$$
does not exist for every $l\ge 1$ with $\phi(l)<(\dim\mathbb{T}^2)^2$.
\end{enumerate}
\end{exa}

Claim (a) is straightforward. 

The eigenvalues of $\hat S$ and $\hat T$ are the primitive roots of
unity of order 4 and 3 respectively. Therefore,
$$
\hbox{Spec}(S^l)\cap \hbox{Spec}(T^l)=\emptyset
$$
unless $l$ is divisible by 12. Since $\phi(l)\ge 4$ for all $l\ge 12$,
this implies (b). 

To prove (c), we take $x_0\in\mathbb{T}^2$ such that the points
$$
x_0, Sx_0, S^2 x_0, S^3 x_0, Tx_0, T^2 x_0
$$
are distinct and a neighborhood $U\subset \mathbb{T}^2$ of $x_0$ such
that
$$
S^nU\cap  T^nU=\emptyset \Leftrightarrow S^n x_0\ne T^nx_0.
$$
Then for $f$ equal to the characteristic function of $U$, we have
$$
\int_{\mathbb{T}^2}f(S^{n}x)f(T^{n}x)\,
dm(x)=\left\{\begin{tabular}{cl}
$m(U)$, & if 12 divides $n$,\\
0, & otherwise.
\end{tabular} \right.
$$
Since $\phi(l)<4$ implies that $l<12$, this proves (c).

\begin{exa}[cf. Corollary \ref{th:two}]\label{ex:two}
There exist (i) infinite-dimensional, (ii) disconnected, compact abelian
group $X$ and epimorphisms $S$ and $T$ of $X$ such that $\{S,T\}$ is not
mixing, but there is no proper closed subgroup $Y$ of $X$ such that
for some $l\ge 1$,  $Y$ is $\{S^l,T^l\}$-invariant 
and $S^l|_{X/Y}=T^l|_{X/Y}$.
\end{exa}

We utilize an example constructed by D. Berend in \cite{ber} for a
different purpose. Let
$$
X=\prod_{n\in\mathbb{Z}} Y
$$
for a compact abelian group $Y$
(with appropriate choice of $Y$, $X$ can be made infinite-dimensional
or disconnected). Note that for every $\chi\in \hat X$, there exists a
finite $D_\chi\subset \mathbb{Z}$ and $\chi_n\in\hat Y-\{0\}$, $n\in D_\chi$,
such that
$$
\chi\left((y_n)_{n\in \mathbb{Z}}\right)=\prod_{n\in D_\chi}
\chi_n(y_n).
$$
Consider the following permutations of $\mathbb{Z}$:
\begin{align*}
\sigma(n)&=n+1,\\
\pi(n)&=\left\{\begin{tabular}{ll} $n$ &for $l_{2k}\le |n|\le l_{2k+1}$,\\
$-n$ &for $l_{2k+1}\le |n|\le l_{2k+2}$,\end{tabular}\right.\\
\tau&=\pi^{-1}\sigma\pi,
\end{align*} 
where $\{l_i\}_{i\ge 1}$ is an increasing sequence of integers such
that $l_0=0$ and 
\begin{equation}\label{eq:l_i}
l_{i+1}/l_i\to\infty\quad\hbox{ as $i\to\infty$.}
\end{equation}
Permutations $\sigma$ and $\tau$ define automorphisms $S$ and $T$ of
$X$ which act on $X$ by permuting coordinates.

First, we observe that $\{S,T\}$ is not mixing. In fact, for 
$$
B=\left\{(y_n)_{n\in\mathbb{Z}}:\, x_0\in A\right\}
$$
where $A$ is a measurable subset of $Y$, we have
\begin{equation}\label{eq:no_lim}
m(S^{-n}B\cap T^{-n}B)=\left\{\begin{tabular}{ll} $m(B)$ & for
    $\pi(n)=n$,\\
$m(B)^2$ & for $\pi(n)\ne n$.\end{tabular}\right.
\end{equation}

Suppose that there exists a proper closed
subgroup $Y$ such that for some $l\ge 1$, $Y$ is $\{S^l,T^l\}$-invariant
and $S^l|_{X/Y}=T^l|_{X/Y}$. This is equivalent to existence of a proper
subgroup $\Gamma$ of $\hat X$ such that $\Gamma$ is $\{\hat S^l,\hat T^l\}$-invariant
and $\hat S^l|_{\Gamma}=\hat T^l|_{\Gamma}$. Consider the map
$$
\hat X\to \{D\subset \mathbb{Z}:\, |D|<\infty\}: \chi\mapsto D_\chi
$$
Since this map is $\left<\sigma,\tau\right>$-equivariant, the image of
$\Gamma$ is a set $\Delta$ consisting of finite subsets of $\mathbb{Z}$ such that
$\Delta$ is $\{\sigma^l,\tau^l\}$-invariant and $\sigma^l(D)=\tau^l(D)$ for every $D\in\Delta$.
It follows that for every $k\ge 1$ and $D\in \Delta$, we have
\begin{equation}\label{eq:D}
\sigma^l\pi\sigma^{kl}(D)=\pi\sigma^{(k+1)l}(D).
\end{equation}
Take $d\in D$. Because of (\ref{eq:l_i}), there exist infinitely many $k_i\ge 1$ such
that 
$$
\pi(\sigma^{k_i l}(d))=\sigma^{k_i l}(d)\quad\hbox{and}\quad \pi(\sigma^{(k_i+1) l}(d))=-\sigma^{(k_i+1) l}(d).
$$
Then by (\ref{eq:D}),
$$
\sigma^{-(k_i+1)l}\pi^{-1}\sigma^l\pi\sigma^{k_il}(d)=-d-2(k_i+1)l\in D.
$$
This contradicts finiteness of $D$. Hence, $\Delta=\{\emptyset\}$ and
$\Gamma=\hat X$ which proves the claim.

\begin{exa}[cf. Corollary \ref{thm:limit}]\label{ex:limit}
There exist (i) infinite-dimensional, (ii) disconnected, compact abelian
group $X$, a Borel subset $B$ of $X$,  and epimorphisms $T$ and $S$ of $X$ such that
for every $l\ge 1$, the limit
$$
\lim_{n\to\infty} m(S^{-ln}B\cap T^{-ln}B)
$$
does not exist.
\end{exa}

Let $S$ and $T$ be as in Example \ref{ex:two}. It follows from
(\ref{eq:l_i}) that for every $l\ge 1$, there exist 
infinitely many $n_1,n_2\ge 1$ such that $\pi(ln_1)=ln_1$ and $\pi(ln_2)\ne ln_2$.
Hence, by formula (\ref{eq:no_lim}), the limit does not exist.

\begin{exa}[cf. Proposition \ref{c:mix}]\label{ex:1}
There exist (i) infinite-dimensional, (ii) disconnected, compact abelian
group $X$ and an infinite subgroup of $\hbox{\rm Aut}(X)$ such that
the action of $\Gamma$ on $X$ is not mixing and every element of infinite
order is ergodic.
\end{exa}
 
Take $$X=\prod_{n\ge 1} Y\quad\hbox{ for a compact
  abelian group $Y$}$$
(choosing $Y$ appropriately, one can make $X$
either disconnected or infinite dimensional).  Take $\Gamma$ to be the group of finitary
permutations of the components of $X$. It is a torsion group which is
not mixing.

To give a less trivial example, consider
\begin{align*}
\Gamma&=\mathbb{Z}\ltimes V\quad\hbox{ with
  $V=\{\pm 1\}^\mathbb{Z}$,}\\
V_0&=\{(v_i)\in V:\, v_i=1\quad\hbox{for $i\ge 1$}\},\\
X&=\prod_{\Gamma/V_0} Y\quad\hbox{ for a compact abelian group $Y$.}
\end{align*}
The group $\Gamma$ acts on $X$ permuting coordinates, and since $V_0$ does
not contain nontrivial normal subgroup, $\Gamma$ embeds in $\hbox{Aut}(X)$.
Every element of infinite order in $\Gamma$ is mixing, but because
$V_0$ is infinite, the action of $\Gamma$ is not mixing.

\begin{exa}[cf. Proposition \ref{c:mix}]\label{ex:2}
There exists a semigroup $\Gamma$ of epimorphisms of the torus $\mathbb{T}^d$
which is not mixing, but its every finitely generated subsemigroup is
mixing.
\end{exa}

Consider
$$
\Gamma=\left<2\cdot\hbox{SL}(d,\mathbb{Z})\right>.
$$
Since
$$
\left(\begin{tabular}{ll} 2 & $2n$\\ 0 & 2\end{tabular}\right)\in \Gamma
$$
for every $n$, it is not mixing. On the other hand, if $\Gamma_0$ is a finitely generated subsemigroup
and
$$
\gamma_i=2^{n_i}\delta_i\in \Gamma_0,\quad \delta_i\in \hbox{SL}(d,\mathbb{Z}),
$$
such that $\gamma_i\to\infty$, then $n_i\to\infty$ as well. This
implies that for every $\chi\in \hat{\mathbb{T}}^d-\{0\}=\mathbb{Z}^d-\{0\}$,
$$
\hat \gamma_i\chi\to\infty.
$$
Hence, the action of $\Gamma_0$ on $X$ is mixing.

\begin{exa}[cf. Corollary \ref{p:bhat}]\label{ex:mix2}
There exists a free nonabelian semigroup $\Gamma$ of epimorphisms of the
torus $\mathbb{T}^d$ which is mixing of all orders.
\end{exa}

Take $\alpha,\beta\in\hbox{SL}(d,\mathbb{Z})$ that generate a free
group and let $\Gamma$ be the semigroup generated by $2\alpha$ and $2\beta$.
It was shown above that $\Gamma$ is mixing.
Suppose that $\Gamma$ is mixing of order $s-1$, but not mixing of order $s$. Then there
exist $x_i\in \mathbb{Z}^d-\{0\}$ and $\gamma_i^{(n)}\in \Gamma$,
$i=1,\ldots,s$, such that $\gamma_j^{(n)}\ldots\gamma_i^{(n)}\to
\infty$ for $1<i\le j\le s$ and
$$
\hat \gamma_1^{(n)} x_1+(\hat \gamma_1^{(n)}\hat \gamma_2^{(n)}) x_2+\cdots +(\hat \gamma_1^{(n)}\ldots\hat \gamma_s^{(n)})x_s=0.
$$
Since $\gamma_2^{(n)}\to\infty$,
$$
\hat \gamma_2^{(n)}=2^{k_n}\delta_n \quad \hbox{for $k_n\to\infty$ and
  $\delta_n\in
\hbox{SL}(d,\mathbb{Z})$.}
$$
It follows that $2^{k_n}$ divides $x_1$ and $x_1=0$. This gives a contradiction.
Hence, $\Gamma$ is mixing of all orders.

\begin{exa}[cf. Corollary \ref{p:bhat} and Lemma ~\ref{c:solv}]\label{ex:mix3}
For $d\ge 4$, there exists a free nonabelian semigroup of automorphisms of the
torus $\mathbb{T}^d$ which generates a solvable group of degree 2 and is mixing of all orders.
\end{exa}

Write $d=d_1+d_2$ with $d_1,d_2\ge 2$, take hyperbolic matrices
$A\in\hbox{SL}(d_1,\mathbb{Z})$, $B\in\hbox{SL}(d_2,\mathbb{Z})$, and consider the semigroup $\Gamma$
generated by
$$
\left(\begin{tabular}{cc} $A$ & $C$\\ 0 & $B$\end{tabular}\right),
\quad C\in \hbox{M}(d_1\times d_2,\mathbb{Z}).
$$
Suppose that there
exist $x_i=(u_i,v_i)\in \mathbb{Z}^d$ and $\gamma_i^{(n)}\in \Gamma$,
$i=1,\ldots,s$, such that $\gamma_j^{(n)}\ldots\gamma_i^{(n)}\to
\infty$ for $1<i\le j\le s$ and
$$
\hat \gamma_1^{(n)} x_1+(\hat \gamma_1^{(n)}\hat \gamma_2^{(n)}) x_2+\cdots +(\hat \gamma_1^{(n)}\ldots\hat \gamma_s^{(n)})x_s=0.
$$
Then
$$
\hat \gamma_1^{(n)}\ldots\hat \gamma_i^{(n)}=
\left(\begin{tabular}{cc} $^t A^{k_i(n)}$ & 0\\ * & $^t B^{l_i(n)}$\end{tabular}\right)
$$
with 
\begin{align*}
k_i(n)\to\infty,\;\;k_{i+1}(n)-k_i(n)\to\infty,\\
l_i(n)\to\infty,\;\;l_{i+1}(n)-l_i(n)\to\infty
\end{align*}
as $n\to\infty$. We have
$$
\sum_{i=1}^s {}^t A^{k_i(n)}u_i=0,
$$
and since $A$ is hyperbolic, $u_i=0$ for $i=1,\ldots, s$. Then
$$
\sum_{i=1}^s {}^t B^{l_i(n)}v_i=0,
$$
and it follows that $x_i=0$ for $i=1,\ldots,s$. This shows that
$\Gamma$ is mixing of all orders. Since matrices $A$ and $B$ are
hyperbolic, the linear map
$$
C\mapsto ACB^{-1},\quad C\in \hbox{M}(d_1\times d_2,\mathbb{Z}),
$$
has eigenvalues $\lambda$ with $|\lambda|\ne 1$. Hence, by
\cite[Theorem~4.17]{ros},
$\Gamma$ contains a free nonabelian semigroup.

\begin{exa}[cf. Proposition \ref{l:irr}]\label{ex:not_irr}
The action of 
$$
\Gamma=\left\{\left(
\begin{matrix}
1 &  *\\
0 & *
\end{matrix}
\right)\right\}\subset \hbox{\rm SL}(d,\mathbb{Z})
$$
on the torus $\mathbb{T}^d$ is ergodic, but not strongly irreducible
and not hereditarily ergodic.
\end{exa}

This is straightforward to check using Proposition \ref{p:ergodic}.

\begin{exa}[cf. Theorem \ref{th:solv1} and Corollary \ref{c:irr}]\label{ex:not_erg}
There exists $\Gamma\subset \hbox{\rm Aut}(\mathbb{T}^3)$ such that 
the action of $\Gamma$ on $X$ is strongly irreducible and in
particular hereditarily ergodic, but $\Gamma$ contains no ergodic elements.
\end{exa}

Consider $\Gamma=\hbox{SO}(2,1)\cap \hbox{SL}(3,\mathbb{Z})$. If the
action of $\Gamma$ on $\mathbb{T}^3$ is not strongly irreducible, then
there exists a subgroup $\Lambda$ of finite index in $\Gamma$ and 
a $\Lambda$-invariant subgroup $A$ of $\mathbb{Z}^3$ such that
$\mathbb{Z}^3/A$ is torsion-free. Then $A\otimes \mathbb{Q}$ is a
proper $\Lambda$-invariant subspace. Since the action of
$\hbox{SO}(2,1)$ on $\mathbb{R}^3$ is irreducible, and $\Lambda$ is
Zariski dense in $\hbox{SO}(2,1)$, this gives a contradiction.
Hence, the action of $\Gamma$ is strongly irreducible.

It is also easy to show that $\Gamma$ contains no ergodic elements.
Denote by $B$ the standard bilinear form and suppose that $\gamma\in
\Gamma$ has no roots of unity as eigenvalues. Let $v,w\in\mathbb{C}^3$
be eigenvectors of $\gamma$ with eigenvalues $\lambda,\mu$
respectively. Then $$B(v,v)=B(\gamma v,\gamma v)=\lambda^2 B(v,v)$$ and
it follows that $B(v,v)=0$. Similarly, $B(w,w)=0$. Since $B$ is
nondegenerate, $B(v,w)\ne 0$. Then the computation as above shows that $\lambda\mu=1$.
This implies that $\gamma$ acts trivially on the orthogonal complement
of the subspace $\left<v,w\right>$, which is a contradiction.

\begin{exa}[cf. Theorem \ref{th:solv1}]\label{ex:her}
There exist a compact connected infinite-dimensional abelian group $X$
  and an automorphism $T$ of $X$ which is mixing, but not hereditarily ergodic.
\end{exa}

Let $Y$ be any compact connected abelian group,
\begin{align*}
X=\prod_{n\in\mathbb{Z}} Y,\quad\hbox{and}\quad
T:(y_n)_{n\in\mathbb{Z}}\mapsto (y_{n+1})_{n\in\mathbb{Z}}.
\end{align*}
Then $T$ is mixing, but $T$ acts trivially on the connected subgroup
$$
\{(y_n)_{n\in\mathbb{Z}}:\, y_n\;\hbox{ is constant}\}.
$$
Hence, $T$ is not hereditarily ergodic.

\begin{exa}[cf. Theorem \ref{th:solv1}]\label{ex:last}
There exist an infinite-dimensional compact connected abelian
group $X$ and an abelian subgroup $\Gamma$ of $\hbox{\rm Aut}(X)$ such
that the action of $\Gamma$ on $X$ is hereditarily ergodic,
but the action of every finitely generated subgroup of
$\Gamma$ is not ergodic. In particular, $\Gamma$ contains no mixing elements.
\end{exa}

Take $T\in\hbox{GL}(2,\mathbb{Z})$ with the characteristic polynomial
$x^2-x-1$. Note that $T$ acts ergodically on the torus $\mathbb{T}^2$.
Consider
\begin{align*}
X&=\prod_{n\ge 1} \mathbb{T}^2,\\
\Gamma &=\prod_{n\ge 1} \left<T\right>\quad\hbox{(direct product).}
\end{align*}
Define $T_i\in\Gamma$, $i\ge 1$, by
$$
T_i\cdot (x_n)_{n\ge 1}=(x_1,\ldots,x_{i-1},Tx_i,x_{i+1},\ldots).
$$
The character group of $X$ is
$$
\hat X=\oplus_{n\ge 1} \mathbb{Z}^2.
$$
We claim that any $\Gamma$-invariant subgroup $S$ of $\hat X$ is of the
form $\oplus_{n\ge 1} S_n$ where $S_n$ is a $\Gamma$-invariant
subgroup of $\hat Y$. Indeed, this follows from the identity
$$
(T_i^2-T_i)\cdot (s_n)_{n\ge 1}=(0,\ldots,0,s_i,0,\ldots),\quad (s_n)_{n\ge 1}\in S.
$$
This implies that any closed connected $\Gamma$ invariant subgroup $Y$
of $X$ has the character group of the form
$$
\hat Y=\oplus_{n\ge 1} \mathbb{Z}^2/S_n.
$$
where $S_n$ is a $T$-invariant subgroup of $\mathbb{Z}^2$ such that
$\mathbb{Z}^2/S_n$ is torsion free, i.e., $S_n=0$ or
$S_n=\mathbb{Z}^2$. Since $T$ acts ergodically on $\mathbb{T}^2$,
the set $\mathbb{Z}^2-\{0\}$ contains no finite $\hat T$-orbits. 
This implies that there are no finite $\Gamma$-orbits in $\hat Y-\{0\}$.
Hence, the action of $\Gamma$ is hereditarily ergodic.

It is easy to see that any finitely generated subgroup of $\Gamma$
fixes some nonzero elements in $\hat X$. Hence, by Proposition
\ref{p:ergodic}, such subgroup is not ergodic.


\begin{thebibliography}{100}
\bibitem{ber0} D. Berend, Ergodic semigroups of epimorphisms.  Trans. Amer. Math. Soc.  289  (1985), 393--407.

\bibitem{ber} D. Berend, Joint ergodicity and mixing.  J. Analyse Math.  45  (1985), 255--284.

\bibitem{ber2}  D. Berend, Multiple ergodic theorems.  J. Analyse Math.  50
 (1988), 123--142.


\bibitem{ben} Y. Benoist, Propri\'et\'es asymptotiques des groupes lin\'eaires.  Geom. Funct. Anal.  7  (1997),  no. 1, 1--47.

\bibitem{bg} V. Bergelson and A. Gorodnik, Weakly mixing group
  actions: a brief survey and an example.  Modern dynamical systems
  and applications,  3--25, Cambridge Univ. Press, Cambridge, 2004.
Corrected version: \texttt{http://front.math.ucdavis.edu/math.DS/0505025}.




\bibitem{bh} S. Bhattacharya, Higher order mixing and rigidity of algebraic actions on compact abelian groups.  Israel J. Math.  137  (2003), 211--221.


\bibitem{br_g} E. Breuillard and  T. Gelander,  A topological Tits
  alternative. To appear in Ann. Math.

\bibitem{cass}  J. W. S. Cassels, Local fields. London Mathematical Society Student Texts, 3. Cambridge University Press, Cambridge, 1986.

\bibitem{eiw} M. Einsiedler and T. Ward, 
Asymptotic geometry of non-mixing sequences.
Ergodic Theory Dynam. Systems  23  (2003),  no. 1, 75--85. 

\bibitem{ew} G. Everest and T. Ward, Heights of polynomials and entropy in algebraic dynamics. Universitext. Springer-Verlag London, Ltd., London, 1999.

\bibitem{koch} H. Koch, Number Theory: Algebraic Numbers and
  Functions. AMS, 2000.

\bibitem{lang} S. Lang, Algebra. Springer-Verlag, New York, 2002.

\bibitem{led}  F. Ledrappier, Un champ markovien peut \^etre
  d'entropie nulle et m\'elangeant.
C. R. Acad. Sci. Paris Ser. A  287  (1978), no. 7, A561--A563.



\bibitem{mas} D. W. Masser, Mixing and linear equations over groups in positive characteristic.  Israel J. Math.  142  (2004), 189--204.

\bibitem{mor} D. Morris, Introduction to arithmetic groups, in preparation;
\texttt{http://people.uleth.ca/\~{}dave.morris/lectures/ArithGrps/}.


\bibitem{pr}  V. Platonov and A. Rapinchuk, Algebraic Groups and
  Number Theory, Academic Press, 1993.

\bibitem{roh1} V.~A.~Rokhlin, 
On endomorphisms of compact commutative groups. (Russian) 
Izvestiya Akad. Nauk SSSR. Ser. Mat. 13 (1949), 329--340.

\bibitem{roh2} V.~A.~Rokhlin, Metric properties of endomorphisms of compact commutative groups. 
Amer. Math. Soc. Transl. 64 (1967), 244--252.


\bibitem{ros} J. Rosenblatt, Invariant measures and growth conditions.  Trans. Amer. Math. Soc.  193  (1974), 33--53.

\bibitem{sch0} K. Schmidt,
Mixing automorphisms of compact groups and a theorem by Kurt Mahler.
Pacific J. Math. 137 (1989), no. 2, 371--385.

\bibitem{sch} K.~Schmidt, Dynamical systems of algebraic origin. Progress in Mathematics 128. Birkh\"auser Verlag, Basel, 1995.


\bibitem{sw} K. Schmidt and T. Ward, Mixing automorphisms of compact groups and a theorem of Schlickewei.  Invent. Math.  111  (1993),  no. 1, 69--76.

\bibitem{sp} T. A. Springer, Linear algebraic groups. Second edition. Progress in Mathematics, 9. Birkha\"user Boston, Inc., Boston, MA, 1998.

\bibitem{tits} J. Tits, Free subgroups in linear groups.  J. Algebra  20  (1972), 250--270.
\end{thebibliography}
\end{document}